\documentclass[preprint,12pt]{elsarticle}

\usepackage{amssymb}
\usepackage{amsmath}   
 \usepackage{amsthm}     

\usepackage{subfigure}   

 \newtheorem{theorem}{Theorem}[section]
\newtheorem{lemma}[theorem]{Lemma}
\newtheorem{proposition}[theorem]{Proposition}
\newtheorem{corollary}[theorem]{Corollary}
\newtheorem{remark}[theorem]{Remark}
\theoremstyle{definition}
\newtheorem{definition}[theorem]{Definition}
\newtheorem{example}[theorem]{Example}

\begin{document}

\begin{frontmatter}



\title{An enhanced basis for producing  B\'{e}zier-like curves}


\author{B. Nouri} 
\author{J.  Saeidian \corref{cor1}}
\ead{j.saeidian@khu.ac.ir}

 \address{Faculty of Mathematical Sciences and Computer, Kharazmi University, No. 50,  Taleghani ave., Tehran, Iran.}

\cortext[cor1]{Corresponding author}

\begin{abstract}
This study aims on proposing a new structure for constructing Bernstein-like bases.
 The structure uses an auxiliary function and a shape parameter to construct a new family of bases from any family of blending functions. The new family of bases  inherit almost all algebraic and geometric properties   of   the initial blending functions. 
The corresponding curves have the freedom to travel from the  curve constructed from the initial blending functions
   to the line segment joining the first and last control points. 
The new bases have the monotonicity preservation property and the
 shape of the curve could be adjusted by changing the parameter.
\end{abstract}




\end{frontmatter}


\section{Introduction}\label{sec1}

Parametric curves and surfaces  play a crucial role in computer aided geometric design (CAGD) and they provide with   powerful and reliable  tools for shape design and geometric representation \cite{far1}. In CAGD  and computer graphics,  curves are often described by means of control points in the form 
\begin{equation}
\mathcal{C}\left(  t\right)  =\sum_{i=0}^{n}F_{i}\left(  t\right)
\mathbf{P}_{i},~t\in\left[  a,b\right]  \subset\mathbb{R}\text{,}%
~n\in\mathbb{N}\text{,} \label{Eq:control_point_based}%
\end{equation}
where $\left\{  \mathbf{P}_{i}\right\}  _{i=0}^{n}\subset\mathbb{R}^{\delta}%
,$\ $\left(  \delta\geq2\right)  $ denote control points and the polyline
determined by them is called the \textit{control polygon}. Sufficiently smooth
real functions $\left\{  F_{i}\right\}  _{i=0}^{n}$\ are called
\textit{blending function}s.

 One of the most  important  families of blending functions  are the Bernstein  polynomials which in turn  lead to the well-known classical B\'{e}zier curves \cite{rida}.   They  have a simple structure and are extensively used in many  fields of engineering and technology. One can find its footprints in robotics and  highway or railway design; technologies like animation design and creation of 3D tensor product surface
models rely on them and they are also employed in font design, satellite path planning and image compression \cite{prau,far2}.

However,  B\'{e}zier models  have a major drawback in data visualization, that is  once the control points are fixed, the shape of the classical B\'{e}zier curve cannot be changed. 
This leads to an important question: How can one modify  the shape of the curve without changing  the given data and control polygon?

An extensive amount of research have been devoted for finding suitable and adjustable curves to approximate a given set of data. One idea, which has been studied by many authors, is  to propose new sets of blending functions equipped with  shape parameter(s), either by modifying the Bernstein polynomials or introducing completely new bases.

We can consider two main approaches, one uses non-polynomial basis functions and the other employs  polynomial bases with shape parameters. Among the non-polynomial techniques, the rational  B\'{e}zier models are the most important, where some  weight factors, for control points, is used to adjust the shape of the curve \cite{main}.

Trigonometric  B\'{e}zier models are also categorized in non-polynomial approaches, Han et al. \cite{Han09} studied the cubic trigonometric  B\'{e}zier curves  for the first time and then variants of trigonometric bases and its hybrid cases have been proposed and employed by scientists \cite{Han15, Zhu, Maqsood2, Bibi2, Ammad}.

C-curves were  proposed by Pottmann \cite{Pott} and then studied by Zhang \cite{Zhang}, they are extensions of cubic curves and  use trigonometric functions in their structure.

They have attracted  so much attention \cite{ Shen} and  have been extended to   H-B\'{e}zier curves \cite{Pott2}, which are known as the  hyperbolic version  of C-curves.

Kov\'{a}cs and V\'{a}rady \cite{Cova17, Cova18} constructed the so called "proximity curves", which are a sequence of curves approaching to the control polygon, but unfortunately      
the basis functions  have  square root expressions so they are quite complex. 
  In  \cite{Juhasz}, Juhasz  presented a  family of   rational B\'{e}zier curves which produce a smooth transition between a B\'{e}zier curve and its control polygon.     

 Another  approach   would be  to construct  polynomial basis
functions with shape parameters. These bases often use the same structure as  classical Bernstein basis functions  along with one, two or even more parameters which play the role of shape controller  \cite{Hu,Yan,L.Yan,Qin1, li, juhasz2018bezier}, so the name Bernstein-like are generally used for them.  

In curves equipped  with shape parameter(s), by changing the parameter(s), the curve
either approaches to the control polygon \cite{Cova17, Cova18, Juhasz} or moves away from the
control polygon \cite{Yan, L.Yan, li}.      In this work, we aim to present a general structure which enables any B\'{e}zier-like curve to provide a
smooth transition  between the curve itself  and the straight line segment joining its first and last control points.

Starting from any blending  functions and  employing an auxiliary function,   we propose a set of  basis functions equipped with a shape parameter. These functions and the curves obtained by them have the most common features of the initial system.  Thanks to  
the  parameter, the shape of the curve can be  adjusted while the control points are fixed.

The  outline of the paper is  as follows:
In Section 2, the new structure is proposed to build new basis functions and their properties are studied; the corresponding parametric curves and their geometric properties are discussed in  Section 3.  Section 4  analyzes some examples and illustrates the geometric effects of the  new structure. Section 5 is devoted to study the monotone interpolation problem and verify the applicability of the new structure in handling this problem.  Finally, Section 6 highlights the main contribution of the study and states some questions for future works.




\section{Basis construction}
Let us consider a system
of sufficiently smooth functions $ \mathcal{F}=\left \{ F_{n,i}:\left [ 0,1 \right ]\rightarrow \mathbb{R} \right \}_{i=0}^{n} $, which have the following   properties
\begin{itemize}
\item [\rm{(a)}] Non-negativity: $ F_{n,i}(t) \geq 0, ~ i = 0, 1, 2,\ldots , n, ~~    \forall t\in[0,1] $.
\item [\rm{(b)}] Partition of unity: $ \sum_{i=0}^{n}F_{n,i}(t)=1,~~  \forall t\in[0,1] $.
\item [\rm{(c)}] Symmetry: $F_{n,i}(t)=F_{n,n-i}(1-t),~  i = 0, 1, 2, \ldots , n, ~~  \forall t\in[0,1] $. 
\item [\rm{(d)}] Endpoint interpolation:
\begin{eqnarray*}
F_{n,i}(0)= 
\begin{cases} 
1, & i=0; 
\\ 
0, & i=1,\cdots,n,
\end{cases}
\hspace{0.7cm}
F_{n,i}(1)= 
\begin{cases} 
1, & i=n; 
\\ 
0, & i=0,\cdots, n-1.
\end{cases}
\label{rr1}
\end{eqnarray*}
\item [\rm{(e)}]Endpoint tangency:
\begin{eqnarray*}
{F}'_{n,i}(0)= 
\begin{cases} 
-m,& i=0; 
\\ 
m,& i=1; 
\\ 
0, & i=2,\cdots,n,
\end{cases}
\hspace{0.7cm}
{F}'_{n,i}(1)= 
\begin{cases} 
m,& i=n; 
\\ 
-m,& i=n-1; 
\\ 
0, & i=0,\cdots, n-2,
\end{cases}
\label{rr1}
\end{eqnarray*}
where $ m$ is a real value dependent on problem conditions.
\end{itemize}
\begin{definition}
Based on  $ \mathcal{F}=\left \{ F_{n,i} \right \}_{i=0}^{n} $, we introduce a new basis  $\mathcal{T}:= \left \{ T_{n,i} \right \}_{i=0}^{n} $, on $[0,1]$  as follows:
\begin{eqnarray}
T_{n,0}\left ( t \right )&=&\left ( 1-\sigma  \right )\left (1-\varphi \left ( t \right )  \right )+\sigma F_{n,0}\left ( t \right ),  \nonumber\\
T_{n,i}\left ( t \right )&=&\sigma F_{n,i}(t), ~~~~~~~~~~~~~~~~~~~~~~~~~~~~~~ i=1,\cdots, n-1,  \label{eq1}\\
T_{n,n}\left ( t \right )&=&\left ( 1-\sigma  \right )\varphi \left ( t \right )  +\sigma F_{n,n}\left ( t \right ),  \nonumber
\end{eqnarray}
where $\sigma \in\left [ 0,1 \right ]$ is a parameter and   the function $ \varphi:\left [ 0,1 \right ]\rightarrow \left [ 0,1 \right ] $ plays the role of an \textit{"auxiliary function"}  which requires the following properties
\begin{equation} \label{fi}
\begin{array}{cll}
 (i)    & \varphi \left ( 0 \right )=0, & ~ \varphi \left ( 1 \right )=1,\\
(ii)   & \varphi \left ( t \right )+ \varphi \left ( 1-t \right )=1, & \\
(iii)    &\frac{d}{dt}\varphi \left ( t \right )\vert _{t=0}=0, & \frac{d}{dt}\varphi \left ( t \right ) \vert _{t=1}=0. 
\end{array}
\end{equation}
\end{definition}
The conditions  (\ref{fi}) are essential for the basis  $\mathcal{T}$ to inherit  the  properties of  $\mathcal{F}$.
The cubic polynomial  $ \varphi(t)=3t^{2}-2t^{3}$  and the trigonometric function    $\varphi(t)= \sin^{2}\left ( \frac{\pi }{2}t \right ) $  satisfy the  conditions of an auxiliary function in this context.  However, the list of auxiliary functions is not limited to these two cases  and there exist other options that could play the role of an auxiliary function. A detailed discussion on possible choices for the auxiliary functions  is reported in Appendix A.

\begin{theorem}
The basis functions generated by Eq. (\ref{eq1}) have the following properties,
\begin{itemize}
\item [\rm{(a)}] Non-negativity: $ T_{n,i}(t) \geq 0$ $ (i = 0, 1, 2,\ldots , n)$.
\item [\rm{(b)}] Partition of unity: $ \sum_{i=0}^{n}T_{n,i}(t)=1 $.
\item [\rm{(c)}] Symmetry: $T_{n,i}(t)=T_{n,n-i}(1-t) $ $(i = 0, 1, 2, \ldots , n)$. 
\item [\rm{(d)}] Endpoint interpolation:
\begin{eqnarray*}
T_{n,i}(0)= 
\begin{cases} 
1, & i=0; 
\\ 
0, & i=1,\cdots,n,
\end{cases}
\hspace{0.7cm}
T_{n,i}(1)= 
\begin{cases} 
1, & i=n; 
\\ 
0, & i=0,\cdots, n-1.
\end{cases}
\end{eqnarray*}
\item [\rm{(e)}]Endpoint tangency:
\begin{eqnarray*}
{T}'_{n,i}(0)= 
\begin{cases} 
-\sigma m,,& i=0; 
\\ 
\sigma m,& i=1; 
\\ 
0, & i=2,\cdots,n,
\end{cases}
\hspace{0.7cm}
{T}'_{n,i}(1)= 
\begin{cases} 
\sigma m,& i=n; 
\\ 
-\sigma m,& i=n-1; 
\\ 
0, & i=0,\cdots, n-2.
\end{cases}
\end{eqnarray*}
\item [\rm{(f)}] Linear independence:
If the system $ \mathcal{F}$ is linearly independent then so is the system $ \mathcal{T}$. 
\end{itemize}
\end{theorem}

\begin{proof}
 We verify each case separately:
\begin{itemize}
\item [\rm{(a)}] 
Non-negativity:
This is obvious from   non-negativity of $F_{n,i}(t)$ and the fact that  $0\leq \varphi(t) \leq 1$.
\item [\rm{(b)}] 
 Partition of unity: It is a straightforward result from the partition of unity of the   basis $\mathcal{F}$.
\item [\rm{(c)}] 
Symmetry: 
 We have $  1-\varphi \left (1-t \right )=\varphi \left (t \right ) $, so the desired symmetry is a result of the  symmetry of $ F_{n,i}(t), ~~  i=0,\cdots, n $.
 \item [\rm{(d)}] 
 Endpoint interpolation: From the endpoint values of $ F_{n,i}\left ( t \right ),  ~ i=0,\cdots, n $  and $ \varphi \left ( 0 \right )=0,~\varphi \left ( 1 \right )=1$, we have  the desired values for  $ T_{n,i}\left ( 0 \right )$ and $ T_{n,i}\left ( 1 \right ), ~ i=0,\cdots, n $.
 \item [\rm{(e)}] Endpoint tangency: It is an obvious result   considering $\frac{d}{dt}\varphi \left ( t \right ) \vert_{t=0}=\frac{d}{dt}\varphi \left ( t \right ) \vert_{t=1}=0 $ and the  endpoint tangency properties of $ F_{n,i}\left ( t \right ),  ~ i=0,\cdots, n $.
 
\item [\rm{(f)}]  Linear independence:  For a fixed $n$, suppose
$\displaystyle \sum_{i=0}^{n}c_{i}T_{n,i}(t)=0$, 
we define the  curve
\begin{equation*}
\mathbf{C}(t)=\sum_{i=0}^{n}\binom{x_{i}}{c_{i}}T_{n,i}(t)=\binom{\sum_{i=0}^{n}x_{i}T_{n,i}(t)}{0},
\end{equation*}
i.e., a curve the control points of which are  on the line  $ y=0 $. 
According to the properties  of  $ T_{n,i}(t) $, $ \mathbf{C}(t) $ passes through the first and last control point. So 
we have $$ \binom{x_{0}}{c_{0}}=\binom{x_{0}}{0}, \binom{x_{n}}{c_{n}}=\binom{x_{n}}{0}, $$
which results in
\begin{eqnarray*}
&&c_{0}=c_{n}=0,\\
&& \sum_{i=1}^{n-1}c_{i}T_{n,i}(t)=\sigma \sum_{i=1}^{n-1}c_{i}F_{n,i}(t)=0.\\
\end{eqnarray*}
Now, if the system $\mathcal{F}$ is  linearly independent and $ \sigma \neq 0 $ we have
\begin{equation*}
c_{i}=0, \hspace{0.5cm} i=1,\cdots, n-1,
\end{equation*}
which completes the proof.
\end{itemize}
\end{proof}


\section{The curve and its properties}

\begin{definition}
	Given control points $\left\{  \mathbf{P}_{i}\right\}  _{i=0}%
	^{n}\in\mathbb{R}^{\delta}~\left(  \delta>1\right)  $, by means of the  new basis
	functions $\left\{  T_{n,i}\right\}  _{i=0}^{n}$, a curve is
	defined in the form
	\begin{equation}\label{Eq:curve}
		\mathbf{C}_{n}^{\sigma}\left(  t\right)  =\sum_{i=0}^{n}T_{n,i}\left(
		t\right)  \mathbf{P}_{i},~t\in\left[  0,1\right]  ,%
	\end{equation}
	where $\sigma\in\left[  0,1\right]  $ is a global shape parameter.
\end{definition}
If we fix a parameter value $t\in\left(  0,1\right)  $ and let the shape
parameter $\sigma$ vary in $\left[  0,1\right]  $, we get the curve along which
the curve point $\mathbf{C}_{n}^{\sigma}\left(  t\right)  $ moves when the shape
parameter is altered. We will refer to this curve as the path of the point
$\mathbf{C}_{n}^{\sigma}\left(  t\right) $.

\begin{proposition}
	\label{Prop:convex_combination}[Convex combination] The curve (\ref{Eq:curve})\ is a convex
	combination of the  curve  determined by the blending system $ \left \{ F_{n,i} \right \}_{i=0}^{n}$ and  the control points
	$\left\{  \mathbf{P}_{i}\right\}  _{i=0}^{n}$ 
	 and the straight line segment  joining $\mathbf{P}_{0}$ and $\mathbf{P}_{n}$, namely%
	\[
	\mathbf{C}_{n}^{\sigma}\left(  t\right)  =\left(  1-\sigma\right)  \mathbf{C}%
	_{n}^{0}\left(  t\right)  +\sigma\mathbf{C}_{n}^{1}\left(  t\right)
	,~t\in\left[  0,1\right]  ,~\sigma\in\left[  0,1\right]  \text{.}%
	\]
	
\end{proposition}

\begin{proof}
	According to Eq. (\ref{Eq:curve}) and Eq. (\ref{eq1}), we have%
\begin{eqnarray*}
\mathbf{C}_{n}^{\sigma}\left(  t\right)  &=&\sum_{i=0}^{n}T_{n,i}\left(t \right)  \mathbf{P}_{i}\\
&=&\left(  1-\sigma\right) \left [ \mathbf{P}_{0}\left ( 1-\varphi \left ( t \right )\right ) +\mathbf{P}_{n}\varphi \left ( t \right ) \right ] \\ 
 &&+\sigma \left [ \sum_{i=0}^{n}F_{n,i}\left(t\right)  \mathbf{P}_{i} \right ].
\end{eqnarray*}
\end{proof}

\begin{corollary}
	Curves $\mathbf{C}_{n}^{\sigma}(t),~\sigma\in\left[  0,1\right],  $ provide a
	 transition between the  curve determined by control points
	$\left\{  \mathbf{P}_{j}\right\}  _{j=0}^{n}$ and $ \left \{ F_{n,i} \right \}_{i=0}^{n} $, and the straight line segment
	terminated by $\mathbf{P}_{0}$ and $\mathbf{P}_{n}$. Moreover the Paths of points of the curve (\ref{Eq:curve}) are straight line segments, that facilitates constrained shape modification.
\end{corollary}

\begin{definition}
	[Monotonicity-preserving system \cite{carnicer1994monotonicity}] A system of real
	functions $\left\{  F_{i}\right\}  _{i=0}^{n}$ is monotonicity-preserving
	(resp., strictly monotonicity-preserving), if for any sequence $\beta_{0}
	\leq\beta_{1}\leq\cdots\leq\beta_{n}$ (resp. $\beta_{0}<\beta_{1}<\cdots
	<\beta_{n}$) in $\mathbb{R}$, the function $\sum_{i=0}^{n}\beta_{i}F_{i}$ is
	increasing (resp. strictly increasing).
\end{definition}

\begin{proposition}[Monotonicity preservation] \label{prop1}
	The function system $\mathcal{T}$ is  monotonicity-preserving if  $\mathcal{F}$ is  monotonicity-preserving and    $ \varphi \left ( t \right ) $ is increasing.

\end{proposition}

\begin{proof}
	It suffices to show that for any sequence $\beta_{0}%
	\leq\beta_{1}\leq\cdots\leq\beta_{n}$, the function $\sum_{i=0}^{n}\beta_{i}T_{i}$ is increasing, i.e.,  the inequality
\begin{equation}
		\sum_{i=0}^{n}\beta_{i}\frac{\mathrm{d}}{\mathrm{d}t}{T}_{n,i}(t)  \geq 0\label{Eq:combination}
	\end{equation}
	holds. We have

\begin{eqnarray}
\sum_{i=0}^{n}\beta_{i}\frac{\mathrm{d}}{\mathrm{d}t}{T}_{n,i}(t)=
\left ( 1-\sigma  \right )\left ( \beta_{n}-\beta_{0} \right ) \frac{\mathrm{d}}{\mathrm{d}t} \varphi \left ( t \right )+\sigma \frac{\mathrm{d}}{\mathrm{d}t} \left (  \sum_{i=0}^{n}\beta_{i}{F}_{n,i}\left ( t \right )\right ).
\end{eqnarray}
On the other hand, because $\mathcal{F}$ is  monotonicity preserving and    $ \varphi \left ( t \right ) $ is increasing, we have $ \frac{\mathrm{d}}{\mathrm{d}t} \varphi \left ( t \right )\geq 0 $   and  $ \frac{\mathrm{d}}{\mathrm{d}t} \left (  \sum_{i=0}^{n}\beta_{i}{F}_{n,i}\left ( t \right )\right ) \geq 0 $. Finally, the desired result follows  from the fact that $ \sigma   \in \left [ 0,1 \right ] $.
\end{proof}

Whenever  the control points are in   $\mathbb{R}^2$ plane, then Proposition \ref{prop1} implies that  the   curve  $\mathbf{C}_{n}^{\sigma}$ (Eq. (\ref{Eq:curve})) is  a monotonicity-preserving curve.

\begin{corollary}\label{pro2}
If  $\mathcal{T}$  is a monotonicity-preserving  system of functions and  $ \left \{\mathbf{P}_{i} \right \}_{i=0}^{n}\subseteq \mathbb{R}^{2} $   are  any  monotone data then the curve $ \sum_{i=0}^{n}\mathbf{P}_{i}T _{i} $ is increasing.
\end{corollary}

The forthcoming two diminishing properties are direct consequences of
monotonicity preservation.

\begin{definition}
	[Length diminishing system] Let denote the length of the control polygon and
	that of the curve (\ref{Eq:control_point_based}) by $L_{p}$ and $L_{c}$,
	respectively. A system of functions $\left\{  F_{i}\right\}  _{i=0}^{n}$ is
	length diminishing, if $L_{c}\leq L_{p}$ for any control polygon.
\end{definition}

\begin{corollary}
	[Length diminution] The function system $\mathcal{T}$ is length diminishing if $\mathcal{T}$ is monotonicity-preserving,
	since monotonicity preservation implies length diminution, cf. Theorem 3.5 of
	\cite{Carnicer1999}.
\end{corollary}

\begin{definition}
	[Hodograph diminishing system] The system of functions $\left\{  F_{i}\right\}
	_{i=0}^{n}$ of (\ref{Eq:control_point_based}) is length diminishing, if the
	convex cone generated by the sides $\left\{  \mathbf{P}_{i+1}-\mathbf{P}%
	_{i}\right\}  _{i=0}^{n-1}$ of the control polygon contains the convex cone
	generated by the tangent vectors of the curve (\ref{Eq:control_point_based}).
\end{definition}

\begin{corollary}
	[Hodograph diminution] If $ \mathcal{T} $ is monotonicity-preserving then the function system $\mathcal{T}$ is hodograph
	diminishing, since a system of functions is monotonicity-preserving if and
	only if it is hodograph diminishing, cf. Theorem 2.7 in \cite{Carnicer1997}.
\end{corollary}


\section{Examples}

In this section, we construct some new curves by the structure presented in Eq. (\ref{eq1}) and visually observe their graphical properties.
\begin{example}

The p-B\'{e}zier basis,  introduced in \cite{Cova18}, has the following representation for  $ n=3 $:
\begin{eqnarray*}
&&M_{0}^{\gamma}(t)=\frac{1}{2}+\frac{3}{2} (r_{1}-r_{0}), \hspace{0.4cm}  M_{1}^{\gamma}(t)=\frac{3}{2} (r_{2}-2r_{1}+r_{0}), \\ &&M_{2}^{\gamma}(t)=\frac{3}{2} (r_{3}-2r_{2}+r_{1}), \hspace{0.4cm}  M_{3}^{\gamma}(t)=\frac{1}{2}-\frac{3}{2} (r_{3}-r_{2}),
\end{eqnarray*}
where $  \gamma, t \in [0, 1] $ and
$$
\begin{array}{ll}
 r_{0}^{\gamma}(t)&=t, \\
  r_{1}^{\gamma}(t)&=\sqrt{\left ( 1-\gamma  \right )(\frac{1}{3}-t)^{2}+\gamma \left ( t-\frac{1}{3}+\frac{2}{3}\left ( 1-t \right )^{3} \right )^{2}}, \\
 r_{2}^{\gamma}(t)&=\sqrt{\left ( 1-\gamma  \right )(\frac{2}{3}-t)^{2}+\gamma \left ( t-\frac{1}{3}+\frac{2}{3}\left [2\left ( 1-t \right )^{3} + 3t\left ( 1-t \right )^{2} \right ] \right )^{2}},\\
  r_{3}^{\gamma}(t)&=1-t.
\end{array}$$

Setting  $F_{3,i}=M_i^{\gamma}, ~i=0,1,2,3$, in Eq. (\ref{eq1})  and getting advantage of the auxiliary function   $ \varphi(t)=3t^{2}-2t^{3}$, we come to the basis functions:
 
\begin{eqnarray}
T_{3,0}\left ( t \right )&=&\left ( 1-\sigma  \right )\left (1-3t^{2}+2t^{3}  \right )+\sigma \left ( \frac{1}{2}+\frac{3}{2} (r_{1}-r_{0}) \right ), \nonumber \\
T_{3,1}\left ( t \right )&=& \frac{3}{2} \sigma (r_{2}-2r_{1}+r_{0}), \label{eqstivan} \\
T_{3,2}\left ( t \right )&=&\frac{3}{2} \sigma  (r_{3}-2r_{2}+r_{1}), \nonumber \\
T_{3,3}\left ( t \right )&=&\left ( 1-\sigma  \right )\left (3t^{2}-2t^{3}  \right )+\sigma \left ( \frac{1}{2}-\frac{3}{2} (r_{3}-r_{2}) \right ).  \nonumber 
\end{eqnarray}
The graphical behavior of curves constructed by the  basis functions  defined in Eq. (\ref{eqstivan}) and  the effect of shape parameter  can be observed in Fig. \ref{figstivan}.  Fig.  \ref{figstivan}(a) shows the spectrum of p-curves for different values of $\gamma$, ranging in $[0,1]$. In Fig.  \ref{figstivan}(b) and Fig.  \ref{figstivan}(c)  one can observe the effect of the new structure in modifying the shape of the p-curves.  The new $ \sigma \gamma$-B\'{e}zier-like curves   provide a smooth transition between the p-B\'{e}zier curve and the line segment joining the first and last control points.
Fig.  \ref{figstivan}(c) shows a  smooth travel of curves between the control polygon and the straight line joining first and last control points,  in this special  case the curve successfully  benefits from two  parameters $\gamma$ and $\sigma$ from two different structures.

\begin{figure} 
\centering
\subfigure[$ \sigma=1 $]{
\includegraphics[width=6cm]{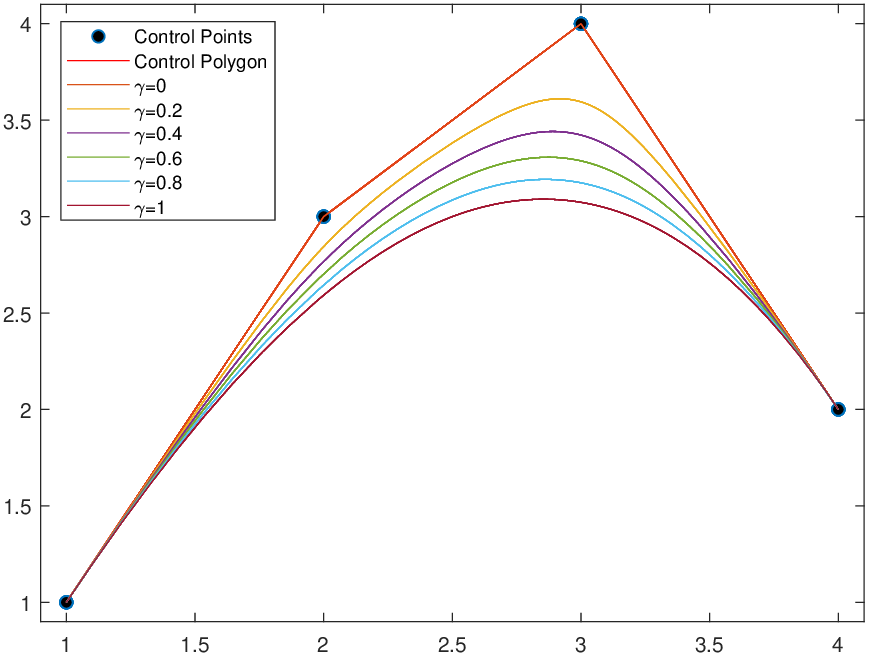}
}
\hspace*{0.1mm}
\subfigure[$ \gamma=1 $ ]{
\includegraphics[width=6cm]{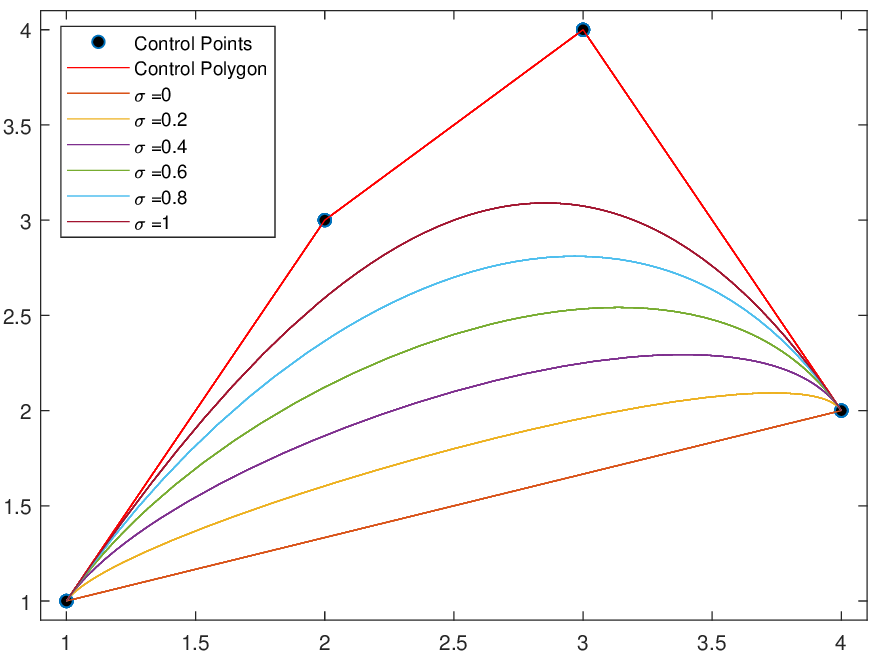}

}
\hspace*{0.1mm}
\subfigure[$ \gamma=0.01 $]{
\includegraphics[width=6cm]{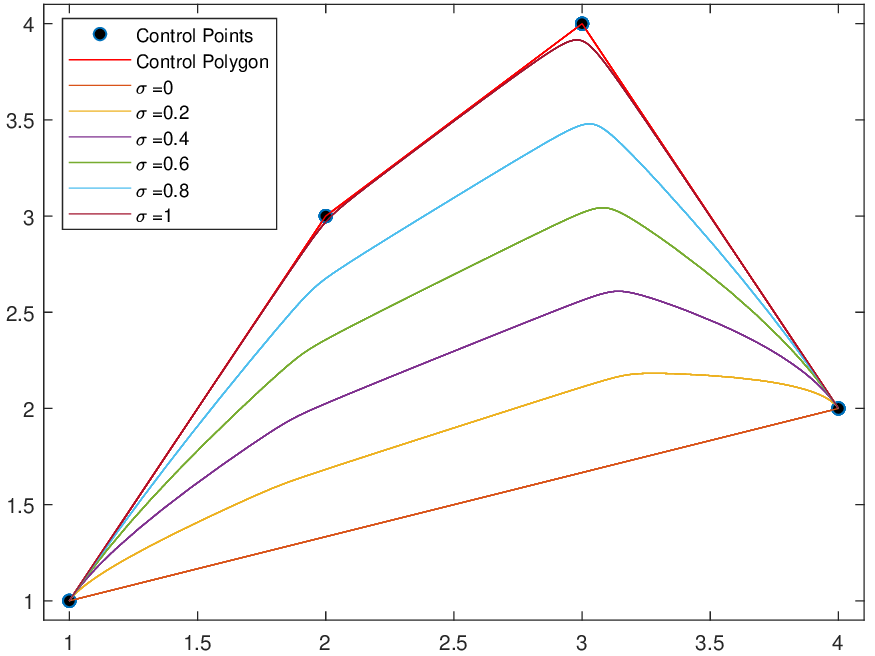}
}

\caption{ (a): The original p-curves with different values of parameter $\gamma$, (b) and (c): the p-B\'
{e}zier curves with $\gamma=1, 0.01$ (resp.) and different values of $\sigma$.} 
\label{figstivan}
\end{figure}

\end{example}


\begin{example}
    Zhu and Han \cite{zhu2014curves}, introduced the   $ \lambda \mu$-Bernstein basis functions,  for shape parameters $ \lambda ,\mu  \in \left [ 0,+\infty  \right ] $, 
\begin{eqnarray}
A_{0}\left ( t; \lambda \right )&=&\left ( 1-t \right )^{3}e^{-\lambda t}, \nonumber \\
A_{1}\left ( t; \lambda \right )&=&\left ( 1-t \right )^{2}\left [ 1+2t-\left ( 1-t \right )e^{-\lambda t} \right ], \label{ezhu} \\
A_{2}\left ( t; \mu \right )&=&\left ( t \right )^{2}\left [ 3-2t-te^{- \mu \left ( 1-t \right )} \right ], \nonumber\\
A_{3}\left ( t; \mu\right )&=&t^{3}e^{-\mu \left ( 1-t \right )}. \nonumber
\end{eqnarray}

The $ \lambda \mu$-curves, constructed by these basis functions, have many basic properties of the  cubic  B\'{e}zier curves, besides by altering shape parameters the curve travels from  the cubic B\'{e}zier curve to the corresponding control polygon. For $ \mu= \lambda = 0$, the $\lambda \mu$-curve reduces to the classical cubic  B\'{e}zier curve.

%
Employing  the structure in Eq. (\ref{eq1}) and getting advantage of the auxiliary function   $$\varphi(t)= \displaystyle \frac{t^{2}}{t^{2}+(1-t)^{2}e^{1-2t}},$$ one can construct a completely new set of blending functions.  

For shape parameters $ \lambda ,\mu  \in \left [ 0,+\infty  \right ] $ and  $ \sigma \in\left [ 0,1 \right ]$, the proposed blending functions $ T_{i}\left ( t; \lambda, \mu, \sigma \right ) $, are defined  for $ t \in\left [ 0,1 \right ] $ as 
\begin{eqnarray}
T_{0}\left ( t; \lambda, \mu, \sigma \right )&=&\left ( 1-\sigma  \right )\left (1-\displaystyle \frac{t^{2}}{t^{2}+(1-t)^{2}e^{1-2t}} \right )   +\sigma \left ( 1-t \right )^{3}e^{-\lambda t}, \nonumber \\
T_{1}\left ( t; \lambda, \mu, \sigma \right )&=&\sigma\left ( 1-t \right )^{2}\left [ 1+2t-\left ( 1-t \right )e^{-\lambda t} \right ],\label{eourzhu}\\
T_{2}\left ( t; \lambda, \mu, \sigma  \right )&=&\sigma t ^{2}\left [ 3-2t-te^{- \mu \left ( 1-t \right )} \right ], \nonumber\\
T_{3}\left ( t; \lambda, \mu, \sigma \right )&=&\left ( 1-\sigma  \right ) \displaystyle \frac{t^{2}}{t^{2}+(1-t)^{2}e^{1-2t}}
+\sigma t^{3}e^{-\mu \left ( 1-t \right )}.  \nonumber 
\end{eqnarray}

 The graphical behavior of the corresponding parametric  curves and  the effect of shape parameter  can be observed in Fig. \ref{ourzhu}. Fig. \ref{ourzhu}(a) and \ref{ourzhu}(b) represent the original curves proposed in \cite{zhu2014curves}, the effect of parameters $\lambda$ and $\mu$ are partially demonstrated in these figures. In Fig. \ref{ourzhu}(c) and \ref{ourzhu}(d)  we  can observe the effect of the new structure, both figures show a smooth transition from the original $\lambda \mu$-curve to the straight line segment joining the last and first control points. 
\begin{figure} 
\centering
\subfigure[$ \mu=0,\sigma=1 $]{
\includegraphics[width=6cm]{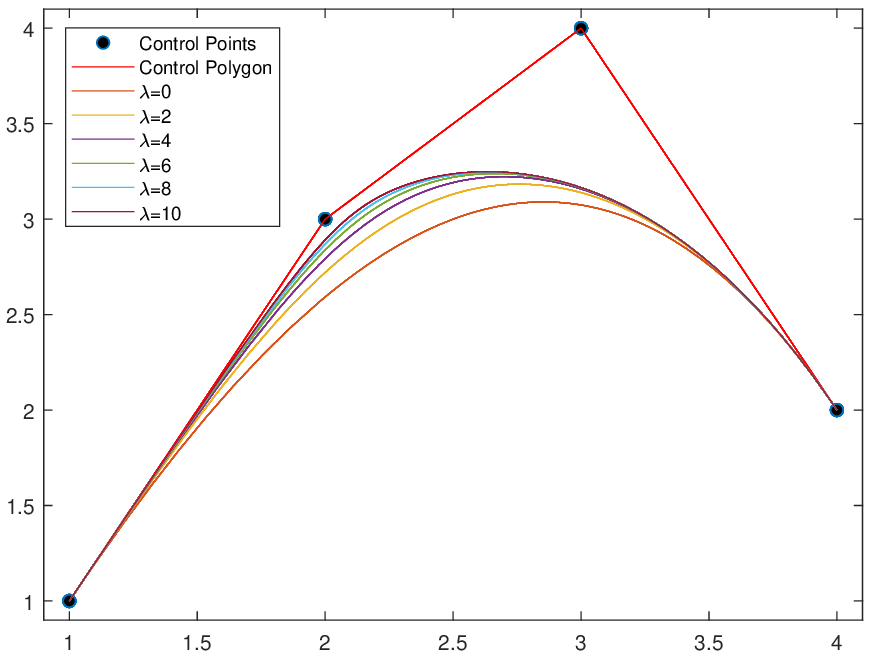}
}
\hspace*{0.1mm}
\subfigure[$\lambda=0,\sigma=1 $ ]{
\includegraphics[width=6cm]{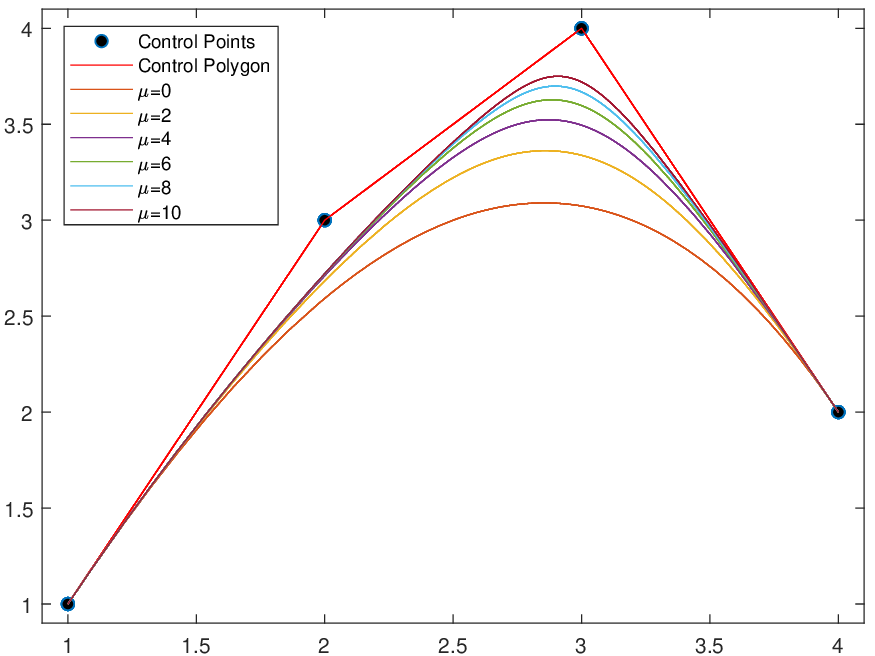}
}

\hspace*{0.1mm}
\subfigure[$ \lambda=0,\mu=0 $ ]{
\includegraphics[width=6cm]{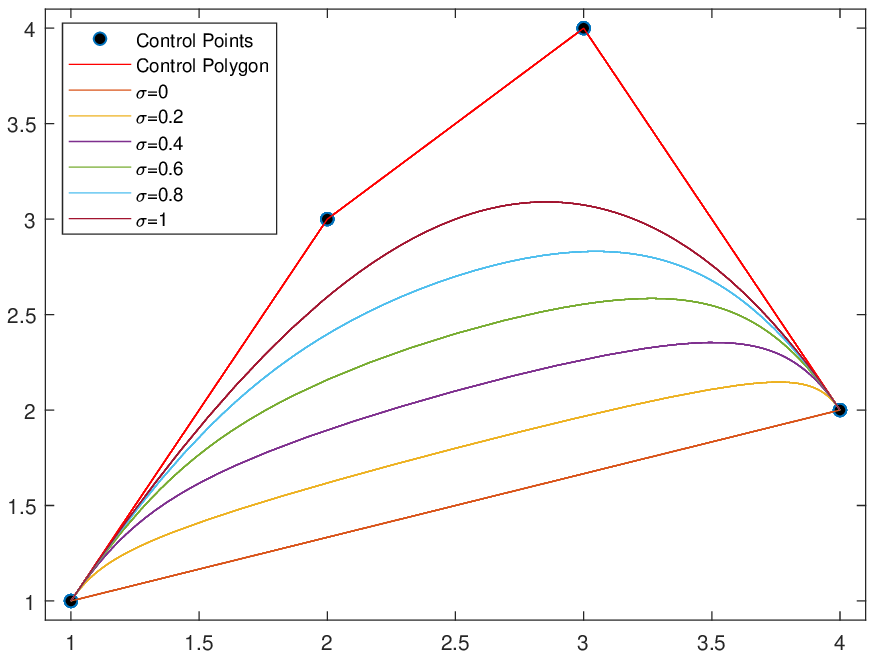}
}
\hspace*{0.1mm}
\subfigure[$ \lambda=10,\mu=10$]{
\includegraphics[width=6cm]{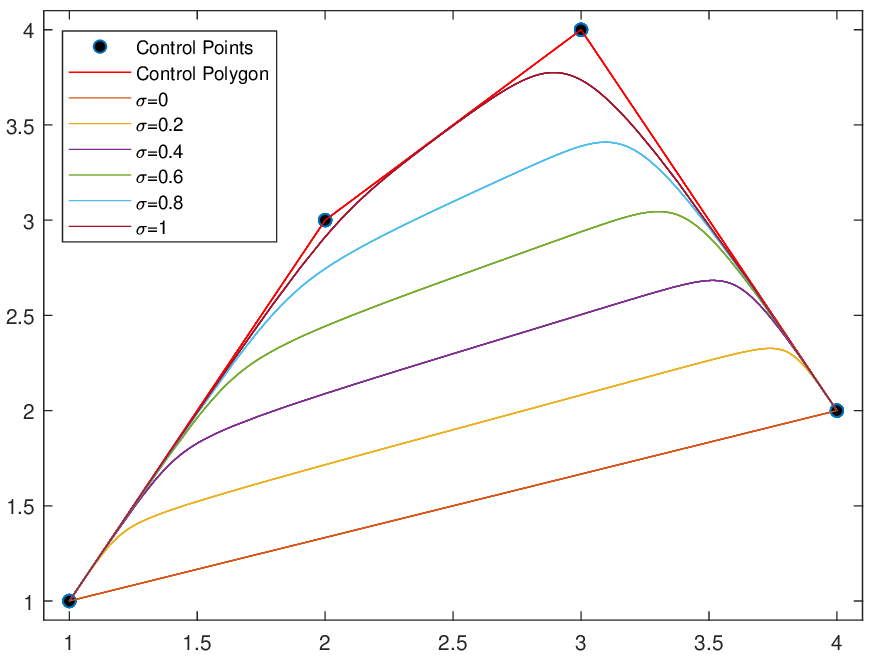}
}
\caption{Curves constructed by basis functions (\ref{eourzhu}), (a) and (b) are the original $\lambda \mu$-curves, and (c) and (d) are the $ \sigma \lambda \mu$-curves for various values of parameters} 
\label{ourzhu}
\end{figure}
\end{example}


\begin{example}\label{ex3}

Let $ \lambda \in \left[ -1,1\right] $, for $ t \in \left[ 0,1\right] $, the Bernstein-like basis functions of order $ 3 $,  introduced in \cite{Yan}, has the following representation:

\begin{eqnarray}
b_{0,3}\left ( t\right )&=& (1-t)^3 (1-2\lambda t +\lambda t^2), \nonumber \\
b_{1,3}\left ( t \right )&=&t(1-t)^2 (3+2\lambda -4 \lambda t +3\lambda t^2), \label{yan2011} \\
b_{2,3}\left ( t \right )&=&t^2 (1-t)(3+\lambda -2 \lambda t +3\lambda t^2) , \nonumber\\
b_{3,3}\left ( t\right )&=&t^3 (1-\lambda  +\lambda t^2). \nonumber
\end{eqnarray}
The  curve constructed by these basis functions have
most of the  properties of the corresponding classical B\'{e}zier curves. When the
shape parameter  increases, the  curve approaches to the control polygon. The new structure enables us to construct  curves  that can approach the line segment connecting the first and last control points. 

Based on the structure proposed in Eq. (\ref{eq1}), and employing the auxiliary functions 
\begin{eqnarray}
(i)&   \varphi_1(t)=&\sin^{2}\left ( \frac{\pi }{2}t \right )\\
(ii)&   \varphi_2(t)=&\displaystyle \frac{t^{2}}{t^{2}+(1-t)^{2}e^{1-2t}}
\end{eqnarray}
one can construct new families of  basis functions.

\textit{Case (i): Trigonometric auxiliary function}

In this case one would have 
\begin{eqnarray}
^1T_{3,0}\left ( t\right )&=& \left ( 1-\sigma  \right )\cos^{2} ( \frac{\pi }{2}t )  +\sigma (1-t)^3 (1-2\lambda t +\lambda t^2), \nonumber \\
^1T_{3,1}\left ( t \right )&=&\sigma t(1-t)^2 (3+2\lambda -4 \lambda t +3\lambda t^2), \label{yan2011NEW} \\
^1T_{3,2}\left ( t \right )&=&\sigma t^2 (1-t)(3+\lambda -2 \lambda t +3\lambda t^2) , \nonumber\\
^1T_{3,3}\left ( t\right )&=&\left ( 1-\sigma  \right )\sin^{2}( \frac{\pi }{2}t )  +\sigma t^3 (1-\lambda  +\lambda t^2). \nonumber
\end{eqnarray}
The formulation has polynomial terms as well as trigonometric ones and it suggests   that  one can increase the effect of periodic terms by setting  smaller values for the shape parameter $\sigma$.
The plots of the corresponding curves are demonstrated in Fig. \ref{YANPIC}, one observes the ability of the new proposed structure to build a spectrum of curves ranging smoothly from the original curve to the line segment joining the extreme control points. The special case $\lambda=1$, which is depicted in Fig. \ref{YANPIC}(d), presents curves with broader range of variations.

\begin{figure} 
\centering
\subfigure[$ \sigma=1 $]{
\includegraphics[width=6cm]{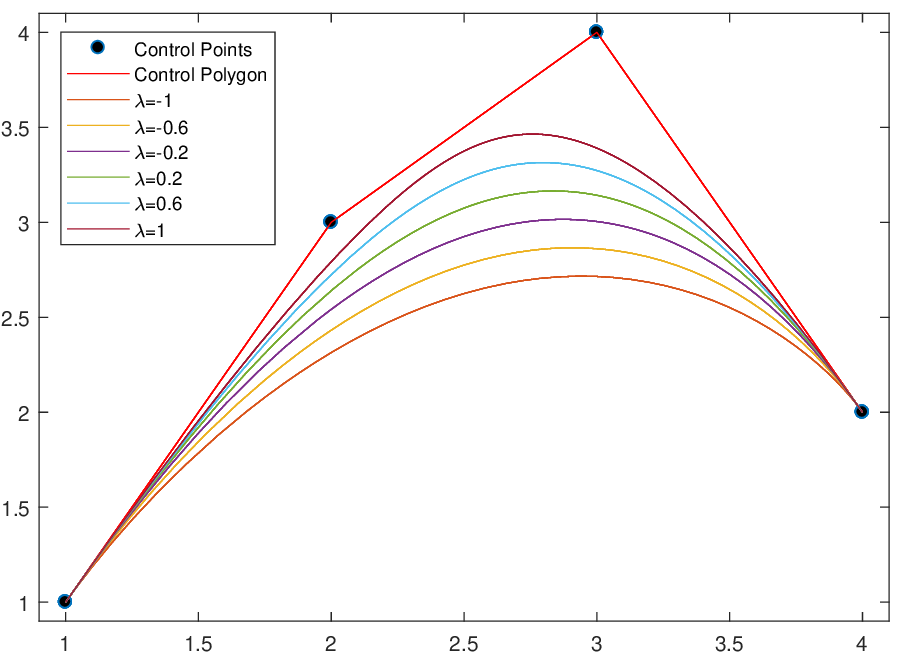}
}
\hspace*{0.1mm}
\subfigure[$ \lambda=0$ ]{
\includegraphics[width=6cm]{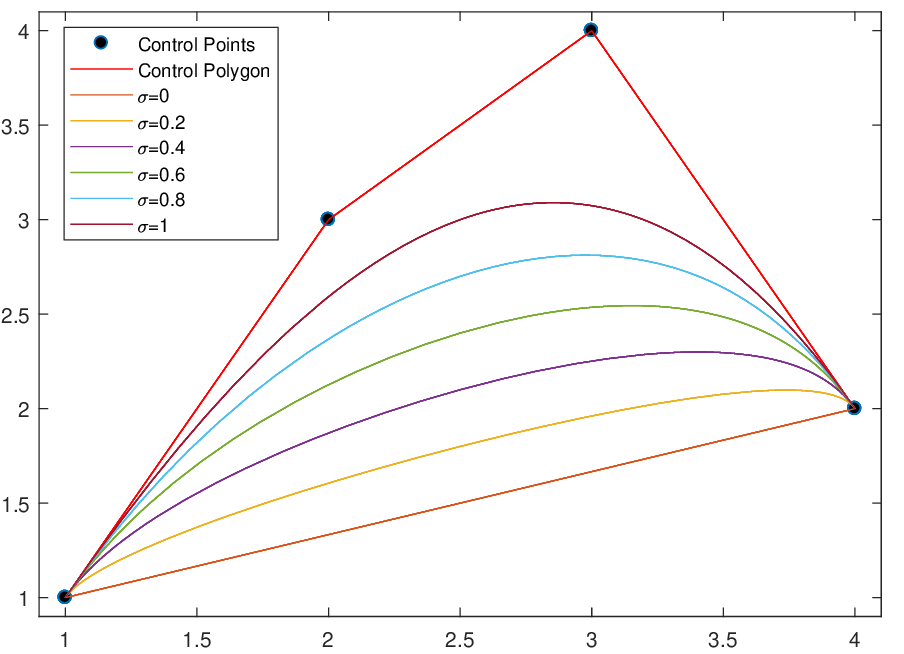}
}
\hspace*{0.1mm}
\subfigure[$ \lambda=-1 $]{
\includegraphics[width=6cm]{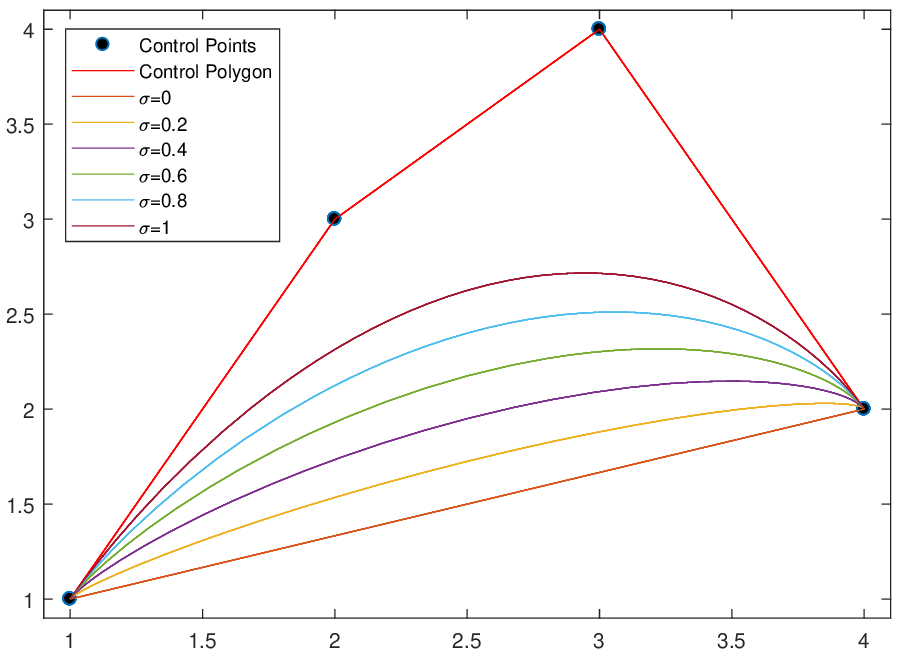}
}
\hspace*{0.1mm}
\subfigure[$ \lambda=1$]{
\includegraphics[width=6cm]{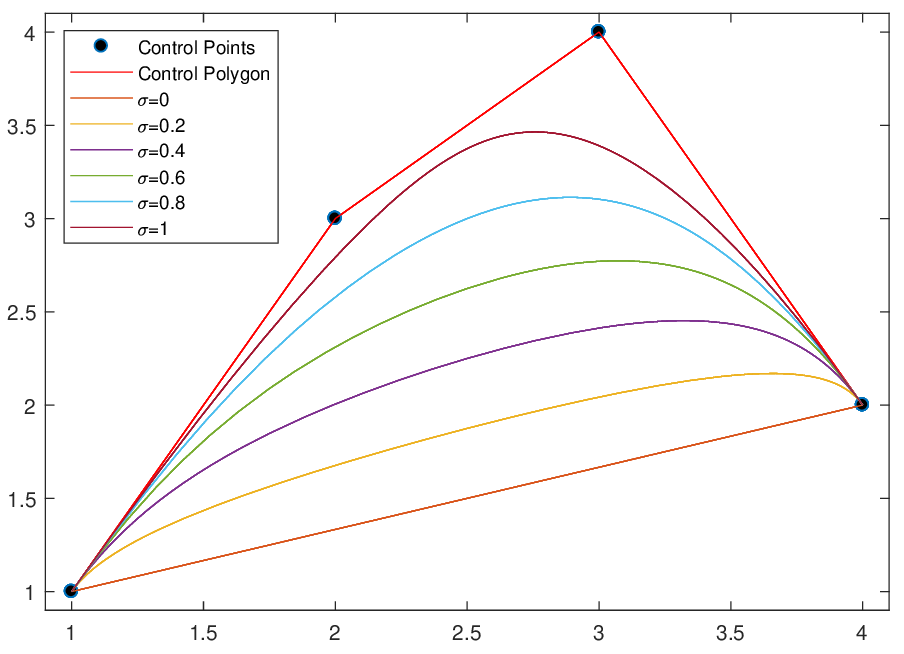}
}

\caption{ Plot of the new curves, according to Example \ref{ex3}-case (i), for different values of parameters} 
\label{YANPIC}
\end{figure}

\textit{Case (ii): Expo-rational auxiliary function}

In this case, the new basis functions would have the formulation
\begin{eqnarray}
^2T_{3,0}\left ( t\right )&=& \left ( 1-\sigma  \right )\left (1- \displaystyle \frac{t^{2}}{t^{2}+(1-t)^{2}e^{1-2t}} \right )   +\sigma (1-t)^3 (1-2\lambda t +\lambda t^2), \nonumber \\
^2T_{3,1}\left ( t \right )&=&\sigma t(1-t)^2 (3+2\lambda -4 \lambda t +3\lambda t^2), \label{yan2011NEW} \\
^2T_{3,2}\left ( t \right )&=&\sigma t^2 (1-t)(3+\lambda -2 \lambda t +3\lambda t^2) , \nonumber\\
^2T_{3,3}\left ( t\right )&=&\left ( 1-\sigma  \right )\left (  \displaystyle \frac{t^{2}}{t^{2}+(1-t)^{2}e^{1-2t}} \right )   +\sigma t^3 (1-\lambda  +\lambda t^2). \nonumber
\end{eqnarray}
This family benefits from having exponential and rational as well as polynomial terms, simultaneously, in its structure.
The corresponding parametric curves are depicted in Fig. \ref{YANPICex}, for  some values of shape parameters.
While Fig. \ref{YANPICex}(a) represents the original curves of \cite{Yan} for different values of $\lambda$, Fig. \ref{YANPICex}(b-d) demonstrate the effect of the new auxiliary function and shape parameter, $\sigma$, in smooth modification of the curves.
\begin{figure} 
\centering
\subfigure[$ \sigma=1 $]{
\includegraphics[width=6cm]{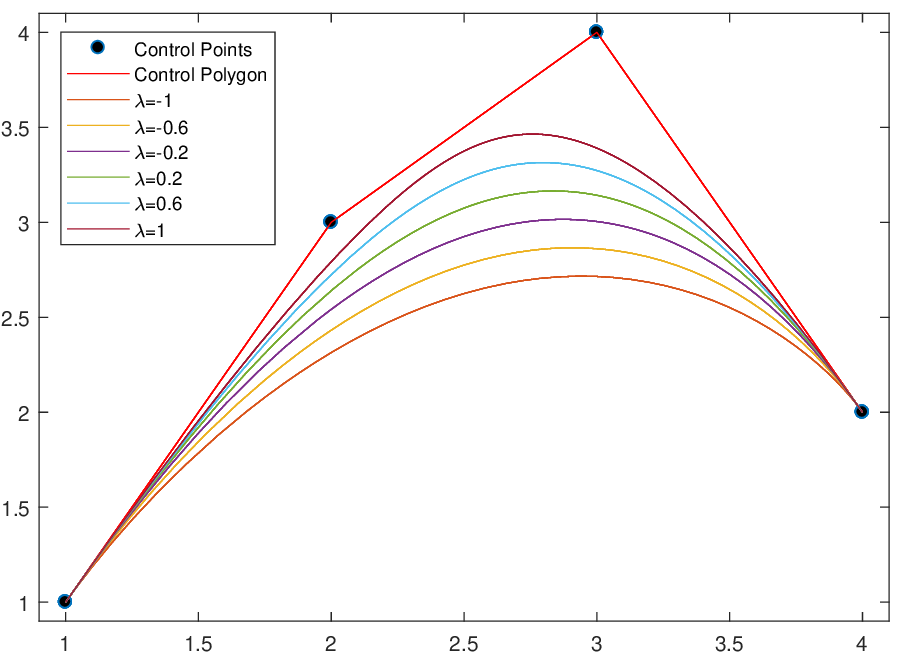}
}
\hspace*{0.1mm}
\subfigure[$ \lambda=0$ ]{
\includegraphics[width=6cm]{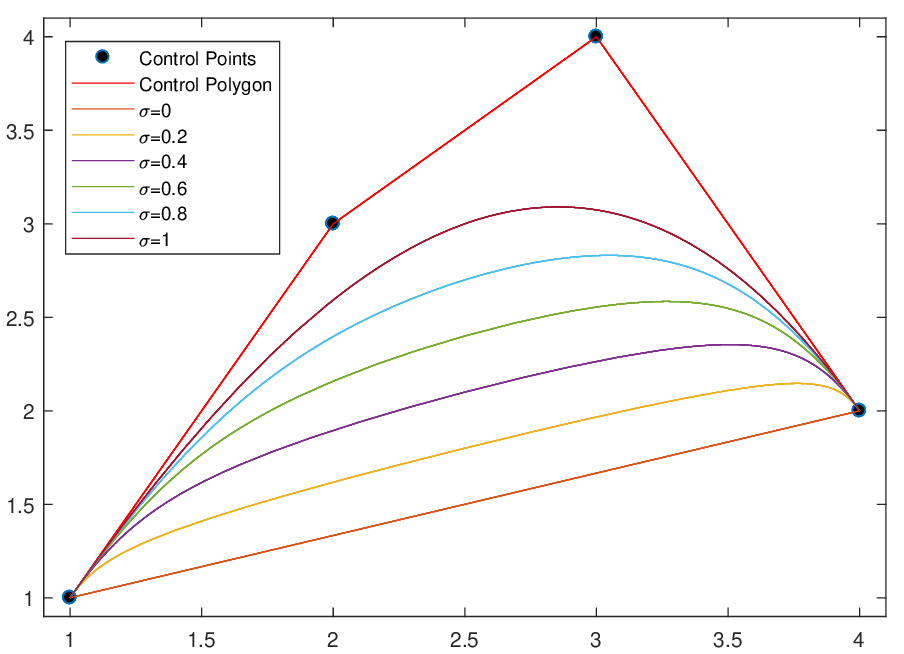}
}
\hspace*{0.1mm}
\subfigure[$ \lambda=-1 $]{
\includegraphics[width=6cm]{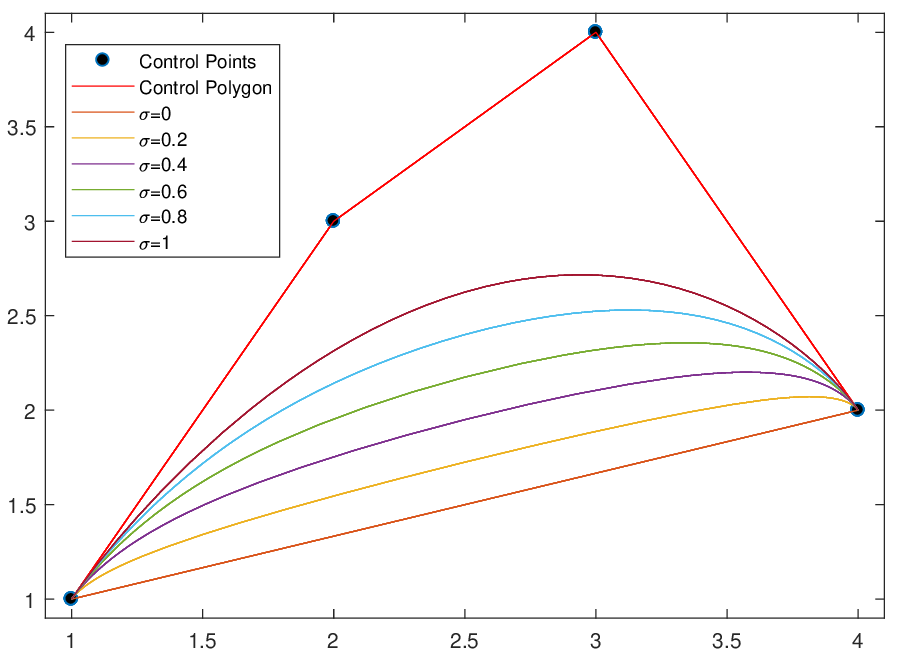}
}
\hspace*{0.1mm}
\subfigure[$ \lambda=1$]{
\includegraphics[width=6cm]{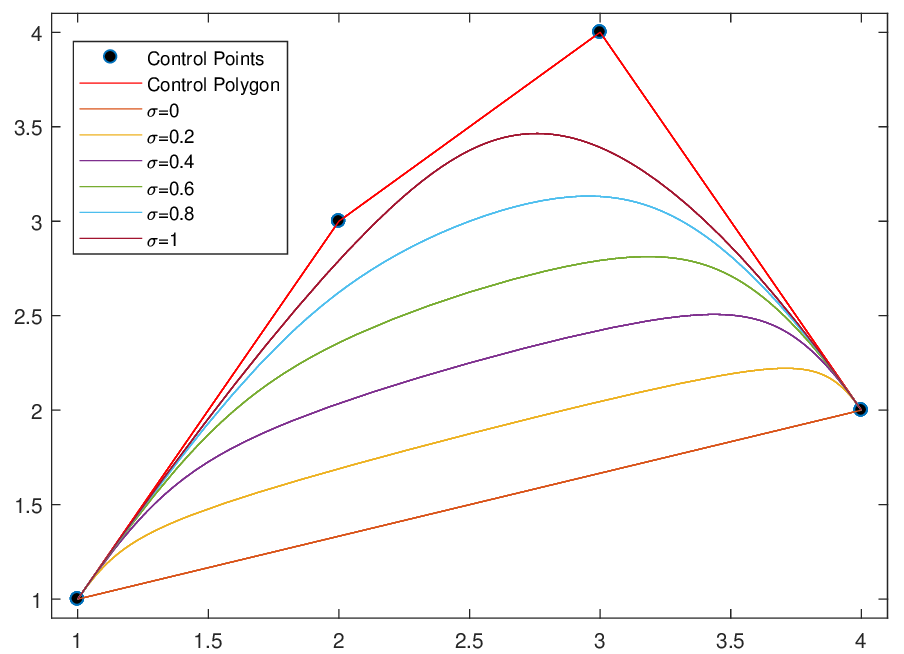}
}
\caption{ Plot of the new curves, according to Example \ref{ex3}-case (ii), for different values of parameters} 
\label{YANPICex}
\end{figure}
\end{example}




\section{ Monotonicity-preserving   interpolation}
Monotonicity-preserving interpolation is a well studied subject in CAD  \cite{car}. 
According to Proposition \ref{prop1}, if there is a monotonicity-preserving system  $\mathcal{F}$,  then getting advantage of an increasing auxiliary function, $\varphi$,  we can construct a new family of monotonicity-preserving  basis functions  $\mathcal{T}$. Now this new system  constructs a family of curves, each of  which   preserves the monotonicity of the control points. Based on this fact,  we employ  the system $\mathcal{T}$ for monotone interpolation and  the shape parameter gives us the freedom to attain $C^1$ or $C^2$ smoothness.

\subsection{$C^1$ Monotone interpolation}
 Let $\left \{ \left ( x_{i},f_{i} \right ) \right \}_{i=0}^{n}$ be a monotone data,  we  present an idea to find a smooth interpolant $p(x)$ on $[x_{0}, x_{n}]$ which  preserves monotonicity of data.

Let   $\mathcal{T}$ be a   monotonicity-preserving system of functions, based on this system, we construct a  piece-wise curve in the form  (\ref{Eq:curve}) as a solution to our monotone interpolation problem.
 For each sub-interval $[x_{i}, x_{i+1}]$ a curve is constructed using the control points
$ \left \{ \left ( x_{i},f_{i} \right ),\left ( h_{i},g_{i} \right ),\left ( t_{i},z_{i} \right ),\left ( x_{i+1},f_{i+1} \right ) \right \} $, where $  h_{i},t_{i},g_{i}$ and $z_{i} $ are unknown values. We need to restrict $ x_{i}< h_{i}< t_{i} <x_{i+1} $ and $  f_{i}\leq g_{i}\leq z_{i} \leq f_{i+1}$  in order to have a monotone curve in each sub-interval  $[x_{i}, x_{i+1}]$. We impose  the piece-wise  curve to be $C^{1}$, so in any two consecutive sub-intervals a continuity condition is needed. These constraints  are illustrated as follows:
 \begin{equation}
\begin{cases} 
t_{i}+h_{i+1}=2x_{i+1} & i=0,\cdots, n-2, \\ 
z_{i}+g_{i+1}=2f_{i+1} & i=0,\cdots, n-2, 
\\
x_{i}< h_{i}< t_{i} <x_{i+1} & i=0,\cdots, n-1,
\\
f_{i}\leq g_{i}\leq z_{i} \leq f_{i+1} & i=0,\cdots, n-1.
\end{cases}\label{15}
\end{equation} 

Whenever  constraints (\ref{15})  have a feasible solution then one obtains a monotone piece-wise
curve which interpolates the given data. 

\begin{theorem}\label{thmc1}

The system of constraints (\ref{15}), always have a feasible solution.

\end{theorem}
\begin{proof}
See Appendix B.
\end{proof}

One feasible solution could be obtained by setting 
\begin{equation}\label{sol1}
\begin{cases}
z_{i}=f_{i+1}, ~~g_i =f_{i},  & i=0, \cdots, n-1,  \\
 t_{i}=x_{i+1}-s,  & i=0, \cdots, n-2,\\
h_{i}=x_{i}+s, & i=1,\cdots, n-1,\\
t_{n-1}=x_{n-1}+2s,~~~~ h_{0}=x_{1}-2s.
\end{cases}
\end{equation}
where  the parameter $ s $ must satisfy:
\begin{equation}\label{fis}
0<s<\min\left \{ \frac{x_{1}-x_{0}}{2},  \frac{x_{2}-x_{1}}{2}, \cdots, \frac{x_{n-1}-x_{n-2}}{2},\frac{x_{n}-x_{n-1}}{2} \right \}.
\end{equation}

\begin{corollary}
For any  monotone set of data $\left \{ \left ( x_{i},f_{i} \right ) \right \}_{i=0}^{n}$,  there exists a piece-wise curve, in the form of Eq. (\ref{Eq:curve}),   which  preserves monotonicity.  It is defined in the sub-interval $[x_{i}, x_{i+1}]$  as:
\begin{equation*}
p_{i}(x)=\binom{x_{i}}{f_{i}}T_{3,0}(t)+\binom{h_{i}}{g_{i}}T_{3,1}(t)+\binom{t_{i}}{z_{i}}T_{3,2}(t)+\binom{x_{i+1}}{f_{i+1}}T_{3,3}(t)
\end{equation*}
where $  h_{i},t_{i},g_{i},z_{i} $  satisfy   constraints (\ref{15}). 
\end{corollary}
\begin{remark}
Constraints (\ref{15}), which  give us the condition of $C^1$ continuity, are satisfied  independent of the value $\sigma$.
  So there is no  need to re-examine the continuity condition $C^1$ by changing $\sigma$. 
\end{remark}
\begin{remark}
The  solution given in Eqs. (\ref{sol1}) and (\ref{fis}) is  just one of the   solutions to the system of constraints  (\ref{15}).   There may be other solutions  and  we will consider those solutions according to  desired constraints. 
\end{remark}

What comes next is an example to show the applicability and reliability of the new structure to handle the monotone interpolation problem.
For demonstration, we get advantage of the  Bernstein basis system $ \mathcal{F}=\left \{ B_{n,i} (t)=\binom{n}{i}t^{i}\left ( 1-t \right )^{n-i} \right \}_{i=0}^{n} $, which  is a monotonicity-preserving system, along with the  the auxiliary  function $\varphi(t) = 3t^2-2t^3$.  This auxiliary function is an increasing one and therefore, according to Corollary \ref{pro2},  the system $\mathcal{T}$,  defined in  equation (\ref{eqmo1}), would be a  monotonicity-preserving system.
\begin{eqnarray}
T_{n,0}\left ( t \right )&=&\left ( 1-\sigma  \right )\left (1-3t^2+2t^3  \right )+\sigma B_{n,0}\left ( t \right ),  \nonumber\\
T_{n,i}\left ( t \right )&=&\sigma B_{n,i}(t), ~~~~~~~~~~~~~~~~~~~~~~~~~~~~~~ i=1,\cdots, n-1,  \label{eqmo1}\\
T_{n,n}\left ( t \right )&=&\left ( 1-\sigma  \right ) \left ( 3t^2-2t^3 \right )  +\sigma B_{n,n}\left ( t \right ).  \nonumber
\end{eqnarray}
\begin{example}\label{mexam7}
 Table \ref{jadvalm7}  represents a  sampling of  the function $f(x)=\displaystyle \frac{1}{1+e^{-x}}$ on $[0, 2]$.
\begin{table}[h!]
\caption{ Data set from Example \ref{mexam7}}\label{jadvalm7}
\centering
\begin{tabular}{|c |c c c c c c c c |}
\hline
$ x $&0&0.292&0.461&0.799&1.172&1.409&1.798&2\\
\hline
$ f(x) $&0.5000&0.572&0.613&0.690&0.763&0.804&0.858&0.881

\\
\hline
\end{tabular}
\end{table}
 We  use the system  defined in (\ref{eqmo1}) to solve the corresponding monotone interpolation problem.
We have solved the corresponding system of constraints (\ref{15}) by setting  $ s=0.05 $, which satisfies (\ref{fis}).
The $C^{1}$ monotonicity-preserving interpolant for  different values of $ \sigma $, as well as the original curve and the error functions are illustrated in Figure \ref{figm7}. The solution is dependent upon the shape parameter $\sigma$, Figure \ref{figm7}(a) illustrates that when the shape parameter decreases the corresponding solution curve gets fewer oscillations. This observation is verified in  Figure \ref{figm7}(b), which presents the error plot.
\begin{figure}[!ht]
\centering
\subfigure[]{
\includegraphics[width=6cm]{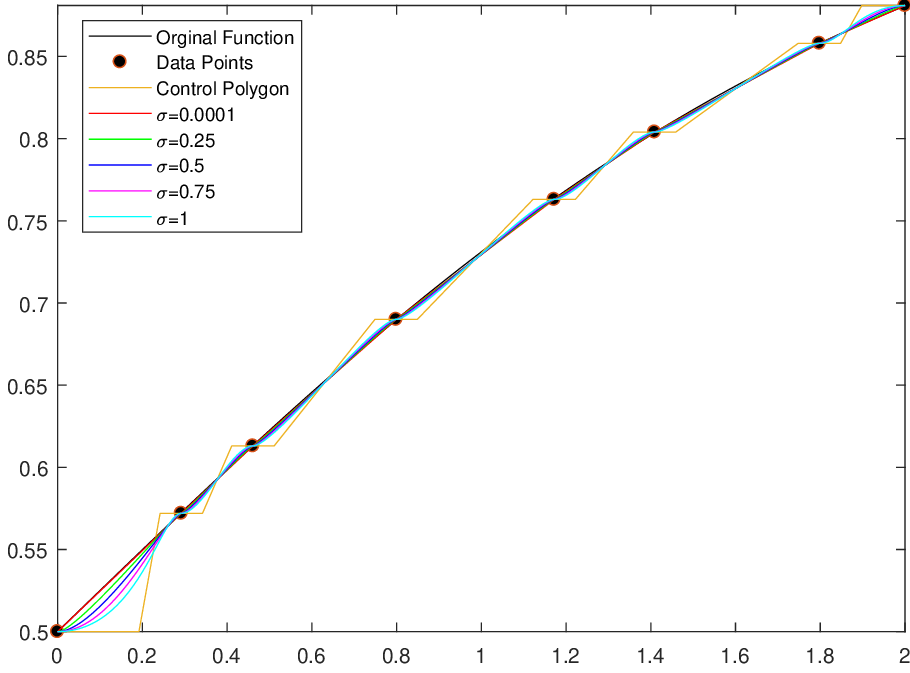}
}
\hspace*{0.1mm}
\subfigure[]{
\includegraphics[width=6cm]{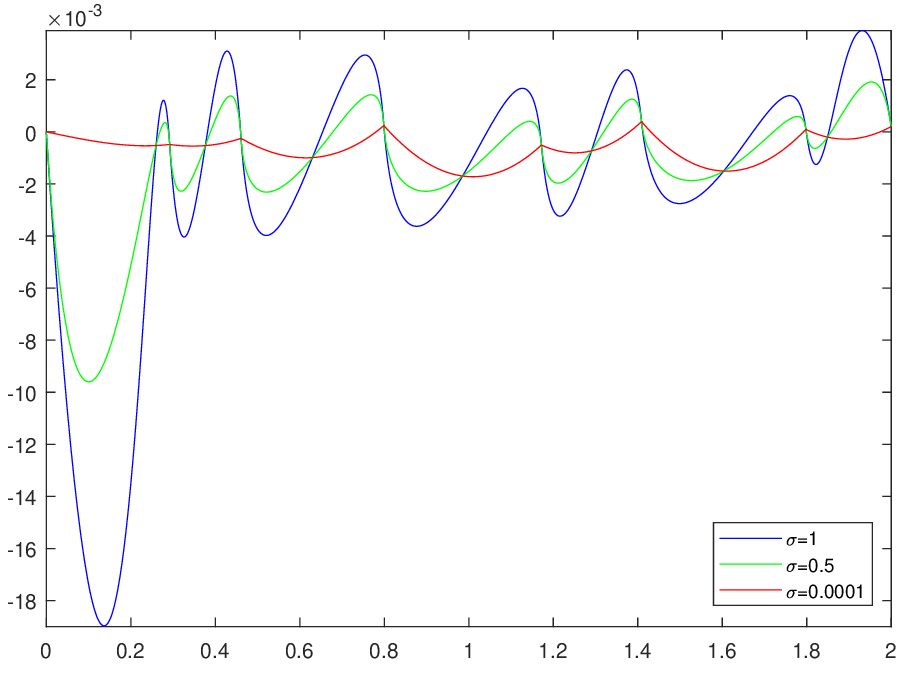}
}

\caption{(a): The $ C^{1} $ monotonicity-preserving interpolant for  different values of $\sigma$, (b): Corresponding error functions. } 
\label{figm7}
\end{figure}
\end{example}

\subsection{$C^2$ monotone interpolation}
It is possible to attain a $C^2$ solution for the monotone interpolation problem. However, this needs some more constraints on the elements of the basis system (\ref{eq1}). 
First of all,   one  needs  the  second derivatives of  $ \mathcal{F}=\left \{ F_{n,i} \right \}_{i=0}^{n} $ to adhere the following endpoint values:
 \begin{eqnarray}
{F}''_{n,i}(0)= 
\begin{cases} 
\omega & i=0;  \\ 
-2\omega & i=1; \\
\omega & i=2;\\
0 & o.w.
\end{cases}
\hspace{0.7cm}
{F}''_{n,i}(1)= 
\begin{cases} 
\omega & i=n;  \\ 
-2 \omega & i=n-1;  \\
\omega & i=n-2;\\
0 & o.w.
\end{cases}
\nonumber
\end{eqnarray}
where $\omega$ is a real value.  
It is also required that the auxiliary function, $\varphi$,  satisfies
 $$ \frac{d^{2}}{dt^{2}}\varphi \left ( t \right ) \vert_{t=0}=    \frac{d^{2}}{dt^{2}}\varphi \left ( t \right ) \vert_{t=1}=0. $$
  In this way,  the second derivatives of  the basis  $ \mathcal{T}=\left \{ T_{n,i} \right \}_{i=0}^{n} $ would read
 \begin{eqnarray}
{T}''_{n,i}(0)= 
\begin{cases} 
\sigma \omega & i=0; 
\\ 
-2\sigma \omega & i=1; 
\\
\sigma \omega & i=2;\\
0 & o.w.
\end{cases}
\hspace{0.7cm}
{T}''_{n,i}(1)= 
\begin{cases} 
\sigma \omega & i=n; 
\\ 
-2\sigma \omega & i=n-1; 
\\
\sigma  \omega & i=n-2;\\
0 & o.w.
\end{cases}
\nonumber
\end{eqnarray}

Now we are ready to attain a $C^2$ interpolant curve which preserves the monotonicity of the data, to do so we will need  extra points as auxiliary points.
At each sub-interval $[x_{i}, x_{i+1}]$, we construct the curve using the control points
\begin{equation*}
\left \{ \left ( x_{i},f_{i} \right ),\left ( h_{i},g_{i} \right ),\left ( t_{i},z_{i} \right ),\left ( w_{i},k_{i} \right ),\left ( d_{i},c_{i} \right ),\left ( x_{i+1},f_{i+1} \right ) \right \},
\end{equation*}
 where $ h_{i},t_{i},w_{i},d_{i},g_{i},z_{i},k_{i},c_{i} $ are unknown values.  The  $C^{2}$ continuity condition is imposed for  two consecutive sub-intervals.
Moreover, we need the  restrictions  $ x_{i}< h_{i}< t_{i}<w_{i}< d_{i} <x_{i+1} $ and $ f_{i}\leq g_{i}\leq z_{i}\leq k_{i}\leq c _{i} \leq f_{i+1}$,  in order  to have a monotone curve in each sub-interval $[x_{i}, x_{i+1}]$,  and so forth  the overall curve would be  monotone. These constraints are illustrated as follows:
\begin{equation}\label{cons2}
\begin{cases} 
d_{i}+h_{i+1}=2x_{i+1}, & i=0,\cdots, n-2, 
\\
c_{i}+g_{i+1}=2f_{i+1}, & i=0,\cdots, n-2, 
\\
2d_{i}-2h_{i+1}=w_{i}-t_{i+1}, & i=0,\cdots, n-2, 
\\
2c_{i}-2g_{i+1}=k_{i}-z_{i+1}, & i=0,\cdots, n-2, 
\\
x_{i}< h_{i}< t_{i}<w_{i}< d_{i} <x_{i+1}, & i=0,\cdots, n-1,
\\
f_{i}\leq g_{i}\leq z_{i}\leq k_{i}\leq c _{i} \leq f_{i+1}, & i=0,\cdots, n-1.
\end{cases}
\end{equation}

\begin{theorem}\label{thmc2}

The system of constraints (\ref{cons2}), always have a feasible solution.

\end{theorem}
\begin{proof}
See Appendix C.
\end{proof}

\begin{corollary}
For a monotone set of data $\left \{ \left ( x_{i},f_{i} \right ) \right \}_{i=0}^{n}$,  there exists a piecewise curve, in the form of Eq. (\ref{Eq:curve}),  that preserves monotonicity.  It is defined in the subinterval $[x_{i}, x_{i+1}]$  as follows:
\begin{equation*}
p_{i}(x)=\binom{x_{i}}{f_{i}}T_{5,0}(t)+\binom{h_{i}}{g_{i}}T_{5,1}(t)+\binom{t_{i}}{z_{i}}T_{5,2}(t)+ \binom{w_{i}}{k_{i}}T_{5,3}(t)+ \binom{d_{i}}{c_{i}}T_{5,4}(t) +\binom{x_{i+1}}{f_{i+1}}T_{5,5}(t)
\end{equation*}
where $ h_{i},t_{i},w_{i},d_{i},g_{i},z_{i},k_{i},c_{i} $ satisfy the  constraints (\ref{cons2}). 
\end{corollary}

\begin{remark}\label{remarkso}
In Appendix C,   the feasibility of the system of constraints   (\ref{cons2}) is verified by introducing a  solution, However
there may be other solutions to (\ref{cons2})   and this gives a degree of freedom in choosing the desired one. 

For example, if the data is strictly monotone, we can present  a solution through   the following relations: 
\begin{itemize}
\item
For $x$-values, we define an auxiliary positive  parameter $\zeta$ and set the values as follows: 
\begin{equation}\label{solx2c2} 
\begin{cases} 
h_{i}=x_{i}+\zeta,   &i=1,\cdots, n-1,
\\
t_{i}=x_{i}+ \frac{5}{2} \zeta,    &i=1,\cdots, n-1,
\\
w_{i}=x_{i+1}- \frac{3}{2} \zeta,    &i=0,\cdots, n-2,
\\
d_{i}=x_{i+1}- \zeta,   &i=0,\cdots, n-2,
\\
h_{0}=x_{i+1}- \frac{5}{2} \zeta,  &
\\
t_{0}=x_{i+1}-2 \zeta,   &
\\
w_{n-1}=x_{n-1}+3 \zeta,   &
\\
d_{n-1}=x_{n-1}+ \frac{7}{2} \zeta,   &
\\
\end{cases}
\end{equation} 
where  $ \zeta $  must satisfy   the following condition
\begin{equation}\label{solx2c2-neq}
0<\zeta<  \frac{1}{4}\min_{i} \left \{ x_{i+1}-x_{i} \right \}.
\end{equation}

\item
Moreover, for $y$-values we define the auxiliary positive parameter $\eta $ and consider the following values
\begin{equation}\label{soly2c2} 
\begin{cases} 
g_{i}=f_{i}+\eta,  &i=1,\cdots, n-1,
\\
z_{i}=f_{i}+ \frac{3}{2} \eta,  & i=1,\cdots, n-1,
\\
k_{i}=f_{i+1}- \frac{5}{2} \eta,   &  i=0,\cdots, n-2,
\\
c_{i}=f_{i+1}- \eta, &   i=0,\cdots, n-2, 
\\
g_{0}=f_{i+1}- \frac{7}{2}\eta, & \\
z_{0}=f_{i+1}-3 \eta, & \\
k_{n-1}=f_{n-1}+2 \eta,  & \\
c_{n-1}=f_{n-1}+ \frac{5}{2} \eta. &  \\
\end{cases}
\end{equation} 

Here  $ \eta $  must satisfy   the following condition
\begin{equation}\label{soly2c2-neq}
0<\eta< \frac{1}{4} \min_{i} \left \{ x_{i+1}-x_{i} \right \}.
\end{equation}
\end{itemize}
\end{remark}

To demonstrate an example we use the Bernstein basis with the auxiliary function $\varphi(t) = 6t^5-15t^4 +10t^{3}$, which results in the following blending system:
\begin{eqnarray}
T_{n,0}\left ( t \right )&=&\left ( 1-\sigma  \right )\left (1-6t^5+15t^4 -10t^{3}  \right )+\sigma B_{n,0}\left ( t \right ),  \nonumber\\
T_{n,i}\left ( t \right )&=&\sigma B_{n,i}(t), ~~~~~~~~~~~~~~~~~~~~~~~~~~~~~~ i=1,\cdots, n-1,  \label{eqmo11}\\
T_{n,n}\left ( t \right )&=&\left ( 1-\sigma  \right )\left (6t^5-15t^4 +10t^{3} \right )  +\sigma B_{n,n}\left ( t \right ),  \nonumber
\end{eqnarray}

\begin{example}\label{exam32}
We use the same data as Example \ref{mexam7}, presented in Table \ref{jadvalm7}  and handle  
the corresponding monotone interpolation problem by employing the basis system (\ref{eqmo11}).

One needs to solve the system of constraints (\ref{cons2}). We present two different sets of solutions:
\begin{itemize}

\item[Solution 1:]
The first proposed solution is constructed   according to the  proposed values  in Appendix C, through Eqs. (\ref{soly1c2})-(\ref{solx1c2-neq}).  We set $s=0.03$ and the corresponding $C^{2}$ monotonicity-preserving interpolant, for different values of $\sigma$, is depicted in Figure  \ref{figm71}(a), and the error function is graphed in Figure  \ref{figm71}(b). 
\begin{figure}[!ht]
\centering
\subfigure[]{
\includegraphics[width=6cm]{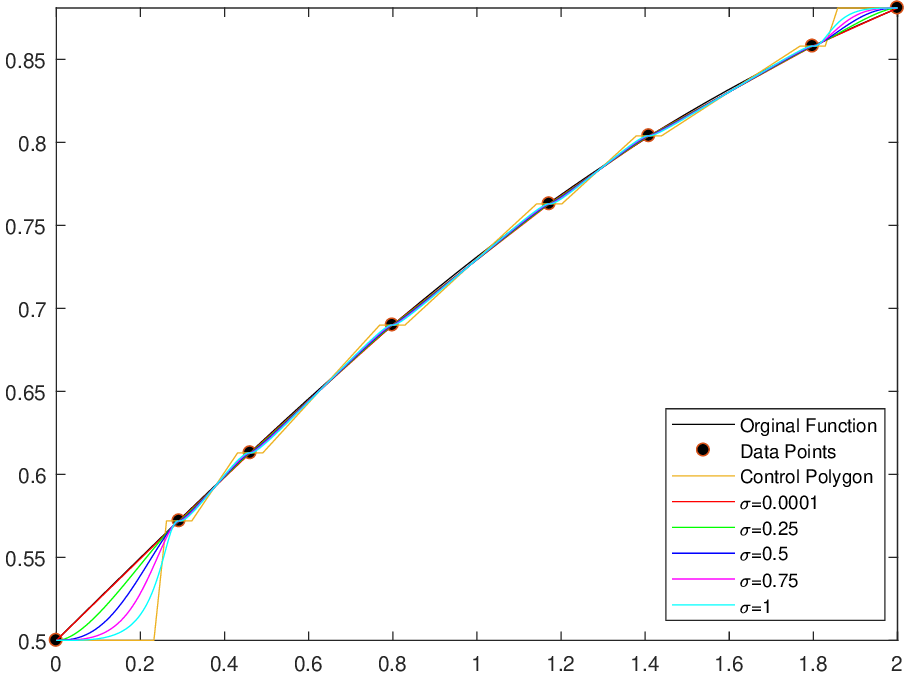}
}
\hspace*{0.1mm}
\subfigure[]{
\includegraphics[width=6cm]{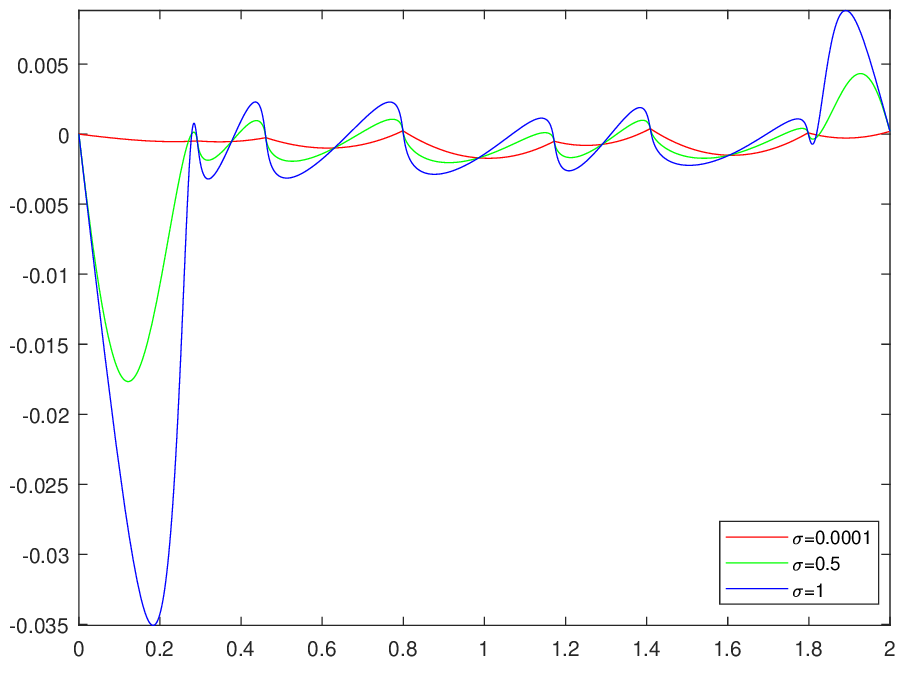}
}

\caption{ Plots  of Solution 1 to Example \ref{exam32}, (a): The $ C^{2} $ monotonicity-preserving interpolant for  different values of $\sigma$, (b): Corresponding error functions. } 
\label{figm71}
\end{figure}

\item[Solution 2:]

Considering  relations Eqs. (\ref{solx2c2})-(\ref{soly2c2-neq}), as a solution to the system (\ref{cons2}), and setting $ \zeta=0.02, \eta=0.003 $,  we come to a series of  
 $C^{2}$ monotonicity-preserving interpolants for  different values of $ \sigma $. The $C^2$ interpolants are pictured in Figure \ref{figm72}-(a)  and the corresponding error functions are plotted in   Figure \ref{figm72}-(b). In comparison with Solution 1, this one 
 shows a significant reduction on the error values.
  Also, we can observe  that Solution 2, competes well with the $C^1$ solution presented in Figure   \ref{figm7}. 

\begin{figure}[!ht]
\centering
\subfigure[]{
\includegraphics[width=5.5cm]{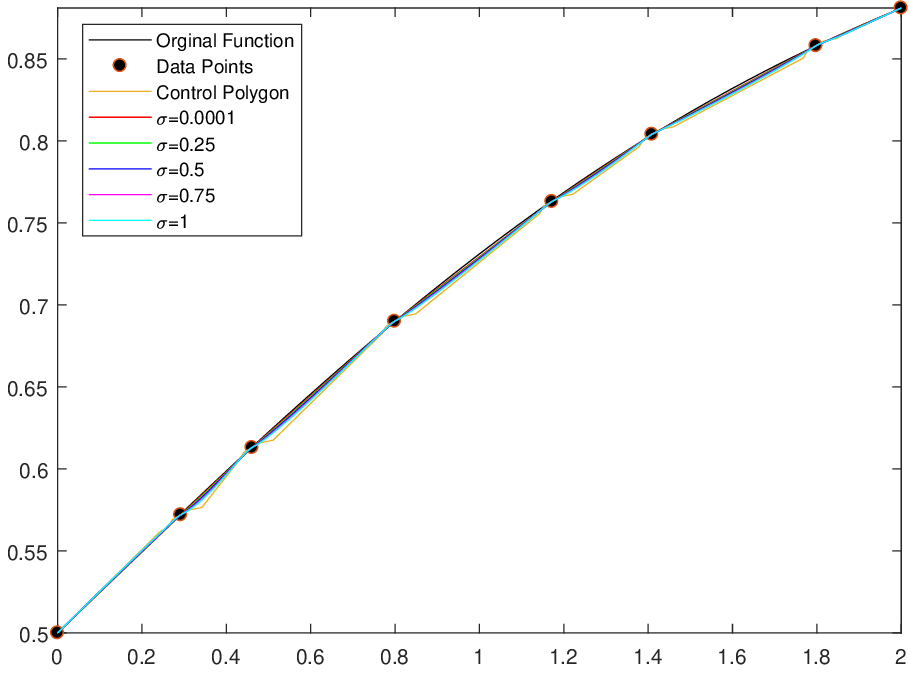}
}
\hspace*{0.1mm}
\subfigure[]{
\includegraphics[width=5.5cm]{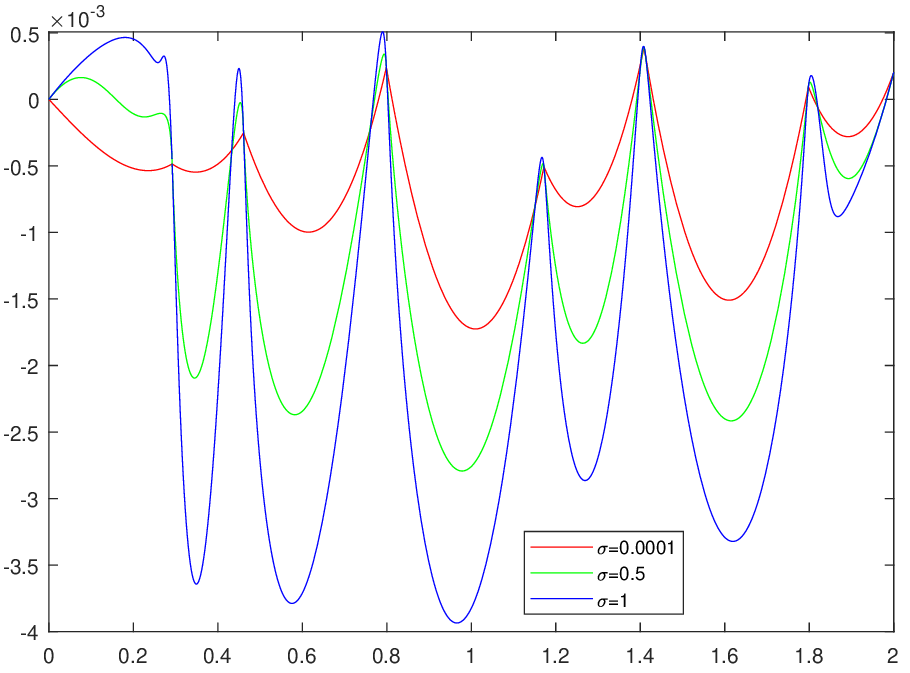}
}

\caption{ Plots  of Solution 2 to Example \ref{exam32}, (a): The $ C^{2} $ monotonicity-preserving interpolant for  different values of $\sigma$, (b): Corresponding error functions. } 
\label{figm72}
\end{figure}

\end{itemize}
\end{example}




\section{Conclusions}

In this study, a general framework  has been proposed to construct Bernstein-like bases. Suppose we have a known set  of blending functions ($\mathcal{F}$), the Bernstein bases, for example;   this new structure uses  an auxiliary function ($\varphi$) and a shape parameter ($\sigma$) to construct a new family of basis functions ($\mathcal{T}$). 
In order to inherit  important algebraic  and geometric  properties of the initial blending functions ($\mathcal{F}$), we imposed some constraints to the auxiliary function.
According to the desired conditions for  the auxiliary function, some suggestions are proposed. This gives a freedom on suitably choosing the  appropriate $\varphi$.  In this way the shape parameter can be used to adjust the shape of the corresponding parametric curves while the control points are fixed.

For any system of  blending functions  as well as the Bernstein basis itself,   this new structure enables us to construct a completely new family of parametric curves equipped with a shape parameter. Under certain conditions the new family is able to provide a tool for solving the monotonicity-preserving interpolation problem.

The quest for finding new auxiliary functions is an open subject. Also, there may be interested readers who want to  extend the study to Bernstein-like operators and apply the new basis for solving PDEs or integral equations.

\section*{Acknowledgments}

We would like to express our  very great appreciation to Professor  Imre Juhasz, head of the Department of Descriptive Geometry, University of Miskolc, Hungary,  for his valuable and constructive suggestions from the 
very first steps of this  work. His willingness to give his time so generously has been greatly appreciated.

\appendix

\section{Some remarks on the auxiliary function $\varphi$}
There exists some other  functions which satisfy the conditions of  the  auxiliary function  $\varphi$, Eq. (\ref{fi}),
here we present some options.

\subsection{Polynomials as $\varphi$ }
We have already introduced the cubic polynomial  $\varphi(t) = 3t^2-2t^3$, as an auxiliary function. 
Here arises a question: {\it{Is there any other polynomial which satisfies  conditions (\ref{fi})?}}

If $n$  is an even number, then there is no polynomial of degree $n$ which satisfies $ \varphi \left ( t \right )+ \varphi \left ( 1-t \right )=1$,  so there would be no even degree polynomial   playing the role of our  auxiliary function.

However, in the case of odd polynomials,  we do have such examples. 
For  $n\geq 3$, an odd integer, consider the function defined by 
$$
\varphi \left ( t \right )=\sum_{j=k+1}^{n}B_{n,j}\left(  t\right) , ~~~~~k=\left\lfloor \frac{n}{2}\right\rfloor,
$$
 where $B_{n,j}\left(  t\right) $ are the  Bernstein polynomials. This is obviously a polynomial of odd degree which satisfies  conditions (\ref{fi}), from now on we call it the \textit{auxiliary polynomial}.

\begin{lemma} 
The auxiliary polynomial is a monotony  increasing  function on $[0,1]$.
\end{lemma}

\begin{corollary}
If $\varphi$ is an auxiliary polynomial and the system $\mathcal{F}$ is monotonicity-preserving then so is  the system $\mathcal{T}$.
\end{corollary}

\begin{remark}
For any even number $n$, the polynomial  
$$\varphi \left ( t \right )  = \frac{1}{2}B_{n,k}\left(  t\right)  +\sum_{j=k+1}^{n}B_{n,j}\left(  t\right),$$
 satisfies the required conditions of an auxiliary function, but after simplifying,  it reduces to a polynomial of degree $n-1$, i.e. 
$$ \frac{1}{2}B_{n,k}\left(  t\right)  +\sum_{j=k+1}^{n}B_{n,j}\left(  t\right)= \sum_{j=k}^{n-1}B_{n-1,j}\left(  t\right).$$

\end{remark}

\subsection{Trigonometric auxiliary functions}

For every odd integer $k$,   function $\varphi_k(t)= \sin^{2}\left ( \frac{k\pi }{2}t \right )$ satisfies the desired conditions of
an auxiliary function.  It is easily seen that for $k=1$,  the trigonometric auxiliary function $\varphi_1(t)= \sin^{2}\left ( \frac{\pi }{2}t \right ),$ is an increasing function on $[0,1]$, which makes it a useful choice  according to Proposition \ref{prop1}.   However, for other values  of $k$, $\varphi_k$ oscillates on  $[0,1]$.


The trigonometric functions reminds the concept of periodicity which models a large amount of  natural phenomenon. With this fact in  mind, suppose one uses a trigonometric auxiliary function with a polynomial basis $\mathcal{F}$ in the framework of Eq. (\ref{eq1}),  this would  result in a basis, $\mathcal{T}$,  which benefits from two important families of basis functions in approximation theory, namely polynomials and trigonometric functions. This  is a verified fact which could be traced back in the literature to the family of trigonometric generalized B-splines \cite{manni2015isogeometric}, just to mention a popular case.

\subsection{Expo-rational auxiliary function}
The function 
$$\varphi(t)= \displaystyle \frac{x^{2}}{x^{2}+(1-x)^{2}e^{1-2x}},$$ 
which is a combination of rational and exponential functions, play the role of an auxiliary function.

\subsection{Pseudo-auxiliary functions}
In the quest for finding suitable auxiliary functions, to be used in the proposed structure (\ref{eq1}), one faces interesting cases which worth mentioning.
There are functions which does not fully satisfy the triple conditions   (\ref{fi}),  but they \textit{nearly} do. One example is the following  function 
 $$\psi(t)= t^{2}\left ( 2\left ( e-2 \right )t+4-e \right )^{t-1},$$ 
which is the product of a quadratic  polynomial and an exponential term. It is a bijection on $[0,1]$ and satisfies conditions   (\ref{fi})-$(i)$ and   (\ref{fi})-$(iii)$, but it fails to satisfy   (\ref{fi})-$(ii)$; however, one has 
 $$0.997 \leq \psi \left ( t \right )+ \psi \left ( 1-t \right )\leq 1,$$
 so, it is very close to be considered  as an auxiliary function.  Accepting $\psi$ as a \textit{pseudo-auxiliary function}, the only property we lose is the partition of unity and thereupon we can not expect the corresponding curve to have the convex hull property.

\section{Proof of Theorem \ref{thmc1}}

We verify the feasibility of  (\ref{15}) by presenting  a solution:
\begin{itemize}
\item 
For $y$-values,  one can set  
\begin{equation*}
z_{i}=f_{i+1},~ g_{i}=f_{i}, \hspace{0.7cm}i=0,\cdots, n-1,
\end{equation*}
which guarantees 
\begin{equation*} 
\begin{cases} 
z_{i}+g_{i+1}=2f_{i+1},     & \hspace{0.7cm} i=0,\cdots, n-2, 
\\
f_{i}\leq g_{i}\leq z_{i} \leq f_{i+1}, & \hspace{0.7cm} i=0,\cdots, n-1.
\end{cases}
\end{equation*}  
\item 
For $x$-values, in order to fulfill the constraints
\begin{equation} \label{xvalues}
\begin{cases}  
t_{i}+h_{i+1}=2x_{i+1} & i=0,\cdots, n-2, 
\\
x_{i}< h_{i}< t_{i} <x_{i+1} & i=0,\cdots, n-1,
\end{cases}
\end{equation} 
we define an artificial variable $s$, which is a real value to be suitably determined.
We set  
\begin{equation*} 
\begin{cases} 
t_{i}=x_{i+1}-s, & i=0,\cdots, n-2, 
\\ 
t_{n-1}=x_{n-1}+2s,
\\
h_{i}=x_{i}+s, & i=1,\cdots, n-1, 
\\
h_{0}=x_{1}-2s.
\end{cases}
\end{equation*} 
In this way the equality constraints in (\ref{xvalues}) hold for any real $s$. On the other hand, from  the inequalities we have
 \begin{eqnarray*}
&&h_{0}=x_{1}-2s < t_{0}=x_{1}-s  \Longrightarrow s>0,\\
&&h_{i}=x_{i}+s < t_{i}=x_{i+1}-s  \Longrightarrow s<\frac{(x_{i+1}-x_{i})}{2}, \hspace{0.5cm}i=1,\cdots, n-2,\\
&&h_{0}=x_{1}-2s\in \left [ x_{0},x_{1} \right ]  \overset{s> 0}{\Longrightarrow}x_{1}-2s\geq x_{0}\Rightarrow s\leq \frac{(x_{1}-x_{0})}{2},\\
&&t_{n-1}=x_{n-1}+2s\in \left [ x_{n-1},x_{n} \right ]\overset{s> 0}{\Longrightarrow}x_{n-1}+2s\leq x_{n}\Rightarrow s\leq \frac{(x_{n}-x_{n-1})}{2}.\\
\end{eqnarray*}
In this way, any $s$ which satisfies 
\begin{equation*}
0<s<\min\left \{ \frac{x_{1}-x_{0}}{2}, \frac{x_{2}-x_{1}}{2},\cdots ,\frac{x_{n-1}-x_{n-2}}{2},\frac{x_{n}-x_{n-1}}{2} \right \},
\end{equation*}
will be a  suitable choice and (\ref{xvalues}) hold.

\end{itemize}

\section{Proof of Theorem \ref{thmc2}}

We prove the theorem by presenting a feasible solution to the system  (\ref{cons2}).
\begin{itemize}
\item 
Set     
\begin{equation}\label{soly1c2}
c_{i}=k_{i}=f_{i+1}~~g_{i}=z_{i}=f_{i}, \hspace{0.7cm}i=0,\cdots, n-1.  
\end{equation}
 In this way, the constraints concerning $y$-values, i.e., 
\begin{equation*} 
\begin{cases} 
c_{i}+g_{i+1}=2f_{i+1} & i=0,\cdots, n-2, 
\\
2c_{i}-2g_{i+1}=k_{i}-z_{i+1} & i=0,\cdots, n-2, 
\\
f_{i}\leq g_{i}\leq z_{i}\leq k_{i}\leq c _{i} \leq f_{i+1} & i=0,\cdots, n-1,
\end{cases}
\end{equation*}  
are  all satisfied.

\item 
For $x$-values, by assigning  
\begin{equation} \label{solx1c2}
\begin{cases}
d_{i}=x_{i+1}+\frac{1}{4}w_{i}-\frac{1}{4}t_{i+1}\\
h_{i+1}=x_{i+1}-\frac{1}{4}w_{i}+\frac{1}{4}t_{i+1}
\end{cases}
\end{equation}
we simply have
\begin{equation*}
\begin{cases} 
d_{i}+h_{i+1}=2x_{i+1} & i=0,\cdots, n-2, 
\\
2d_{i}-2h_{i+1}=w_{i}-t_{i+1} & i=0,\cdots, n-2, 
\end{cases}
\end{equation*} 
so the equality constraints concerning $x$-values in (\ref{cons2}) are fulfilled.
However, yet we  must have
\begin{equation*}
x_{i}< h_{i}< t_{i}<w_{i}< d_{i} <x_{i+1} \hspace{0.5cm} i=0,\cdots, n-1.
\end{equation*}
To impose these constraints, we get advantage of an artificial variable $s$ and set as follows:
\begin{equation} \label{solx1c2-eq}
\begin{cases} 
w_{i}=x_{i+1}-s & i=0,\cdots, n-2, 
\\ 
w_{n-1}=x_{n-1}+2s,
\\
t_{i}=x_{i}+s & i=1,\cdots, n-1, 
\\
t_{0}=x_{1}-2s,
\\
h_{0}=x_{1}-3s,
\\
d_{n-1}=x_{n-1}+3s.
\end{cases}
\end{equation} 
Now according to the inequalities:
\begin{equation*} 
\begin{cases} 
d_{i}<x_{i+1} & i=0,\cdots, n-2, 
\\
h_{i}>x_{i} & i=1,\cdots, n-1, 
\\
w_{i}<d_{i} & i=0,\cdots, n-1, 
\\
h_{i}<t_{i} & i=0,\cdots, n-1, 
\end{cases}
\end{equation*} 
we must have

\begin{eqnarray*}
&&t_{0}=x_{1}-2s < w_{0}=x_{1}-s    \Longrightarrow    s>0,\\
&&h_{0}=x_{1}-3s\in \left [ x_{0},x_{1} \right ]\overset{s> 0}{\Longrightarrow}x_{1}-3s\geq x_{0}\Longrightarrow s\leq \frac{(x_{1}-x_{0})}{3},  \\
&&d_{n-1}=x_{n-1}+3s\in \left [ x_{n-1},x_{n} \right ]\overset{s> 0}{\Longrightarrow}x_{n-1}+3s\geq x_{n}\Longrightarrow s\leq \frac{(x_{n}-x_{n-1})}{3},\\
&&t_{i}=x_{i}+s < w_{i}=x_{i+1}-s \Longrightarrow s<\frac{(x_{i+1}-x_{i})}{2}, \hspace{0.5cm}i=1,\cdots, n-2,\\
\end{eqnarray*}
which naturally leads to 
\begin{equation}\label{solx1c2-neq}
0<s<\min\left \{ \frac{x_{1}-x_{0}}{3}, \frac{x_{2}-x_{1}}{2},\cdots ,\frac{x_{n-1}-x_{n-2}}{2},\frac{x_{n}-x_{n-1}}{3} \right \}.
\end{equation}
Because of $ x_{0}<x_{1}<\cdots<x_{n} $, one always have a feasible  value for  $ s $, which in turn results in a solution to the system:
\begin{equation*} 
\begin{cases}  
d_{i}+h_{i+1}=2x_{i+1} & i=0,\cdots, n-2, 
\\
2d_{i}-2h_{i+1}=w_{i}-t_{i+1} & i=0,\cdots, n-2, 
\\
x_{i}< h_{i}< t_{i}<w_{i}< d_{i} <x_{i+1} & i=0,\cdots, n-1,
\end{cases}
\end{equation*}
 and so  the feasibility of the system of constraints (\ref{cons2}) is verified.

\end{itemize}









\bibliographystyle{elsarticle-num} 
 \bibliography{AddingFunction}
\end{document}